\documentclass[12pt,reqno]{amsart}
\usepackage{amsmath, amssymb, enumitem, float, mathrsfs, amsaddr, graphicx}
\usepackage{tikz-cd}
\usepackage[headings]{fullpage}

\usepackage{nicematrix}

\usepackage{tabularray}
\UseTblrLibrary{booktabs}

\usepackage{booktabs}
\usepackage{changepage}
\usepackage{colortbl}
\usepackage{arydshln}
\usepackage{comment}

\usepackage{hyperref}

\numberwithin{equation}{section}

\theoremstyle{plain}
\newtheorem{theorem}{Theorem}[section]

\newtheorem{proposition}[theorem]{{Proposition}}
\newtheorem{corollary}[theorem]{{Corollary}}

\newtheorem{question}[theorem]{{Question}}

\theoremstyle{definition}

\newtheorem{example}[theorem]{{Example}}

\theoremstyle{remark}
\newtheorem{remark}[theorem]{{Remark}}

\newenvironment{salign}{
    \begin{equation}
    \begin{aligned}
}{
    \end{aligned}
    \end{equation}
    \ignorespacesafterend
}

\newenvironment{graytext}{\color{gray}}{\ignorespacesafterend}

\renewcommand{\pmod}[1]{\ (\mathrm{mod}\ #1)}

\def \A {{\mathbb A}}

\def \E {{\mathbb E}}
\def \F {{\mathbb F}}
\def \G {{\mathbb G}}

\def \P {{\mathbb P}}
\def \Q {{\mathbb Q}}
\def \R {{\mathbb R}}

\def \V {{\mathbb V}}
\def \Z {{\mathbb Z}}

\def\bT  {{\mathbb T}}

\def\cE {{\mathcal E}}

\def\cG {{\mathcal G}}

\def\cN {{\mathcal N}}
\def\cO {{\mathcal O}}
\def\cP {{\mathcal P}}
\def\cS {{\mathcal S}}

\def\cX {{\mathcal X}}
\def\cY {{\mathcal Y}}

\def \sG {{\mathscr G}}

\def \sP {{\mathscr P}}

\def \a {{\mathfrak a}}

\def \c {{\mathfrak c}}

\def \p {{\mathfrak p}}
\def \q {{\mathfrak q}}

\def \w {{\textnormal{\textbf{w}}}}

\newcommand{\GL}{\textnormal{GL}}

\newcommand{\aff}{\textnormal{aff}}

\newcommand{\Aut}{\textnormal{Aut}}

\newcommand{\Gal}{\textnormal{Gal}}
\newcommand{\Hom}{\textnormal{Hom}}
\newcommand{\Ht}{\textnormal{Ht}}

\newcommand{\res}{\textnormal{res}}
\newcommand{\sep}{\textnormal{sep}}

\newcommand{\Sel}{\textnormal{Sel}}

\newcommand{\Val}{\textnormal{Val}}

\newcommand{\mult}{\textnormal{mult}}

\newcommand{\defeq}{\stackrel{\textnormal{def}}{=}}

\begin{document}

\title{Unbounded average Selmer ranks of elliptic curves in  torsion families}

\author{Tristan Phillips}
\email{tristanphillips72@gmail.com}


\begin{abstract}
 	Let $M$ and $N$ be positive integers for which the modular curve $X_1(M,MN)$ has genus $0$, and let $p$ be a prime divisor of $MN$. This article gives asymptotic lower bounds for the average size of the $p$-Selmer group of elliptic curves over a number field, with torsion subgroup $\Z/M\Z\oplus \Z/MN\Z$. In many cases it is shown that this average is unbounded. 
\end{abstract}

\maketitle

\section{Introduction}

Let $K$ be a number field, and let $E$ be an elliptic curve over $K$. 
Throughout this article, we order elliptic curves by \textit{naive height},\footnote{As defined in equation (\ref{eq:naive_height}).} which we denote by $\Ht(E)$.
The average rank of elliptic curves over $K$ is conjectured to be $1/2$. 
A conjecture of Bhargava, Kane, Lenstra, Poonen, and Rains have conjectured that the average size of the $n$-Selmer group of elliptic curves is the sum of the positive divisors of $n$ \cite{PR12, BKLPR15}. Proving this conjecture for arbitrarily large $n$ would imply that the average rank is $1/2$. 
Over the rational numbers $\Q$, Bhargava and Shankar have  computed the average size of the $n$-Selmer group of elliptic curves for each $n\in\{2,3,4,5\}$ \cite{BS15a, BS15b, BS13a, BS13b}. 
Using their results on the average size of the $5$-Selmer group, they were able to conclude that the average rank of elliptic curves over $\Q$ is bounded above by $0.885$. 

In this article, we study families of elliptic curves for which the average size of the $n$-Selmer group can be unbounded. Let $\cY_1(M,N)$ denote the moduli stack of elliptic curves, whose $K$-points parametrize triples $(E,P,C)$, where
\begin{itemize}
\item $E$ is an elliptic curve over $K$,
\item $P\in E[MN](K)$ is an $NM$-torsion point of $E$, and
\item $C\leq E[M](\overline{K})$ is a $\Gal(\overline{K}/K)$-stable cyclic subgroup of order $M$, such that $$E[M](\overline{K})=\langle NP\rangle \oplus C.$$
\end{itemize} 

Let $\cX_1(M,MN)$ denote the smooth compactification of $\cY_1(M,MN)$, and let $X_1(M,MN)$ denote the coarse space of $\cX_1(M,MN)$. 
For $B\in \R_{>0}$, let $\cE_{(M,MN)}(B;K)$ denote the set of elliptic curves $E$ defined over $K$ with height $\Ht(E)\leq B$ and for which there exists $P$ and $C$ with $(E,P,C)\in \cY_1(M,MN)(K)$.

Let $\cO_K$ denote the ring of integers of $K$.
For $m\in \Z_{>1}$, define $\sG_m(K)$ to be the subgroup of $\left(\Z/m\Z\right)^\times$ consisting of elements $a\in \left(\Z/m\Z\right)^\times$, for which there exist infinitely many prime ideals $\p\subset \cO_K$ with norm $N_{K/\Q}(\p)=p$ a rational prime satisfying $p\equiv a\pmod{m}$. 

Let $\Sel_{m}(E)$ denote the $m$-Selmer group of $E$. 
The main result in this article is the following theorem concerning the average size of $\Sel_{m}(E)$ in families.

\begin{theorem}\label{thm:main}
Suppose that $X_1(M,MN)$ has genus zero, and suppose that $p$ is a prime divisor of $MN$. Then there exists an explicit constant $\theta$ (see Table \ref{tab:mean_variance_II}), depending only on $M$, $N$, $K$, and $p$, such that
\begin{equation}\label{eq:main}
	\limsup_{B\to\infty} \frac{1}{\# \cE_{(M,MN)}(B; K)} \sum_{E\in \cE_{(M,MN)}(B;K)} \# \Sel_p(E) \gg \log(B)^{\theta}.
\end{equation}
If additionally
\begin{equation}
(M,N, p, \sG_{MN}(K))\not\in 
\begin{Bmatrix}
(1,5, 5, (\Z/5\Z)^\times), (1,7, 7, \{1,4,7\} \text{ or } (\Z/7\Z)^\times), \\
 (1, 10, 5, (\Z/10\Z)^\times), (2,2, 2, \ast), (2, 6, 3, (\Z/6\Z)^\times), \\ (3,3, 3, \ast), (3,6, 3, \ast), (4,4,2,(\Z/4\Z)^\times), (5,5,5,\ast)
\end{Bmatrix},
\end{equation}
then $\theta$ is positive, implying that the average size of the $p$-Selmer group is unbounded.
\end{theorem}

\begin{remark}
	Equation (\ref{eq:main}) of Theorem \ref{thm:main} can be extended to $d$-Selmer groups for all divisors $d\in \Z_{>1}$ of $M$, and examples of this are given in Section \ref{sec:Composite_case}. 
\end{remark}

Theorem \ref{thm:main} generalizes several results in the literature. In 2014, Klagsbrun and Lemke Oliver proved that the average size of the $2$-Selmer group of elliptic curves over $\Q$ with a $2$-torsion point is unbounded \cite[Corollary 2]{KLO14}. In 2019, Chan, Hanselman, and Li proved that the average size of the $2$-Selmer group of elliptic curves over $\Q$ with $\Z/2\Z\times \Z/8\Z$ torsion  is unbounded \cite[Theorem 6.3]{CHL19}. In a recent preprint of Chan and Verzobio, results similar to Theorem \ref{thm:main} are obtained for many one-parameter families of elliptic curves over $\Q$ \cite{CV25}, including families with prescribed torsion subgroups.\footnote{At present, some of the constants we obtain differ from those in \cite{CV25}. We hope to have this resolved in the near future.}

As in \cite{KLO14, CHL19, CV25}, the strategy is to study the distribution of the ratios of sizes of certain Selmer groups, and prove that they satisfy an analog of the celebrated Erd\H{o}s--Kac theorem \cite{EK40}. The key new ingredient is a systematic approach to counting elliptic curves over arbitrary number fields having a prescribed torsion subgroup and satisfying prescribed local conditions. This is achieved by extending previous work of the author on counting points of bounded height, satisfying various local conditions, on modular curves \cite{Phi24+}.

\subsection*{Contents} 
In Section \ref{sec:preliminaries} we introduce notation and collect relevant facts about number fields, distributions of norms of prime ideals in congruence classes, modular curves, weighted projective stacks, heights, isogenies, and Selmer groups. 
In Section \ref{sec:Selmer_ratios} we study ratios of certain Selmer groups and see how they are related to ratios of certain Tamagawa numbers.
In Section \ref{sec:local_conditions} we explain how results of \cite{Phi24+} can be extended to count elliptic curves over number fields with a prescribed torsion subgroup and prescribed local conditions (Proposition \ref{prop:elliptic_count_general}). 
 In Section \ref{sec:Selmer_distribution} we prove that the logarithms of ratios of sizes of certain Selmer groups are normally distributed (Theorem \ref{thm:selmer_erdos_kac_prime}), and use this to prove our main result (Theorem \ref{thm:main2}). 
In Section \ref{sec:Composite_case} we give a couple of examples extending Theorem \ref{thm:main} to $d$-Selmer groups when $d$ is composite.
Finally, Section \ref{sec:tables} contains tables of relevant data concerning the families of elliptic curves we consider.

\subsection*{Acknowledgments}

Many thanks to John Voight for introducing me to the paper \cite{KLO14} and helping me compute defining equations for some of the modular curves, to Asher Auel for helpful conversations about Galois cohomology, and to Robert Lemke Oliver for helpful conversations about counting elliptic curves and encouragement. This research was supported by National Science Foundation, via grant DMS-2303011.

\section{Preliminaries}\label{sec:preliminaries}

In this section we collect various definitions, notation, and results that will later be used.

\subsection{Number fields} 

Let $K$ be a number field of degree $\deg_{\Q}(K)$, let $\cO_K$ denote the ring of integers of $K$, and let $\overline{K}$ be an algebraic closure of $K$. For an ideal $\a\subseteq \cO_K$ let $N_{K/\Q}(\a)\defeq \#(\cO_K/\a)$ denote the norm of $\a$. For a non-zero prime ideal $\p\subset \cO_K$, let $\F_\p\defeq \cO_K/\p$ denote the residue field at $\p$. 
Let $\Val(K)$ denote the set of places of the number field $K$, and let $\Val_0(K)$ and $\Val_\infty(K)$ denote the set of finite and infinite places of $K$, respectively. 
 
 For each place $v\in \Val(K)$, let $K_v$ denote the completion of $K$ with respect to $v$. If $v\in \Val_0(K)$ is a finite place, let $K_v$ denote the completion of $K$ with respect to $v$, let $\pi_v$ be a uniformizer for $K_v$, let $\cO_{K,v}\defeq \{a\in K_v : v(a)\geq 0\}$ denote the valuation ring, let $\p_v$ denote the prime ideal of $\cO_K$  corresponding to $v$, let $\F_v\defeq \F_{\p_v}$ denote the residue field, and let $q_v\defeq \#\F_\p$ denote the norm of $\p_v$. Let $m_v$ denote Haar measure on $K_v$, normalized so that $m_v(\cO_{K,v})=1$.
 
 For any field $F$, let $\cG_F\defeq \Gal(F^{\sep}/F)$ denote the absolute Galois group of $F$.

\subsection{Distribution of norms of prime ideals in congruence classes}

 Let $m\in\Z_{> 1}$ be an integer greater than $1$, and let $a$ be an integer relatively prime to $m$. Let $\sP(K)$ denote the set of non-zero prime ideals of $K$. Let $\p\subset\cO_K$ be a prime ideal lying above a rational prime $p$. The \textbf{degree} of $\p$ over $p$ is $\deg(\p)\defeq \log_p(\#(\cO_K/\p))$. Let $\sP_1(K)$ denote the set of non-zero prime ideals of $K$ of degree $1$. Let $\sP(B;K)$ (resp., $\sP_1(B;K)$) denote the set of primes $\p$ in $\sP(K)$ (resp., in $\sP_1(K)$) for which $N_{K/\Q}(\p)\leq B$. The natural density of $\sP_1(K)$ in $\sP(K)$ is  $1$, i.e., 
\begin{equation}\label{eq:degree_1_primes}
\lim_{B\to\infty}\frac{\#\sP_1(B;K)}{\#\sP(B;K)}=1
\end{equation}
(see, e.g., \cite[Corollary 1 and Corollary 2 of Proposition 7.17]{Nar04}).
 
 Define the counting function
\begin{equation}
\pi(B;K,a,m)\defeq \#\{\p\subset \sP : N_{K/\Q}(\p)\equiv a\pmod{m},\ N_{K/\Q}(\p)\leq B\}.
\end{equation}
Let $\zeta_m$ be a primitive $m$-th root of unity, and let $\deg_K(K(\zeta_m))$ denote the degree of the extension $K(\zeta_m)/K$.

\begin{proposition}\label{prop:K-Dirichlet}
Define
\begin{equation}
\varepsilon(K,a,m)	\defeq 
\begin{cases}
1 & \text{if for infinitely many } \p\in \sP_1(K)\\ 
&\ \text{ we have } N_{K/\Q}(\p)\equiv a \pmod{m},\\
0 & \text{otherwise.}
\end{cases}
\end{equation}
Then
\begin{equation}
\pi(B;K,a,m)\sim \frac{\varepsilon(K,a,m)}{\deg_K(K(\zeta_m))}\frac{B}{\log(B)}.
\end{equation}
\end{proposition}

\begin{proof}
This follows from \cite[Lemma 7.31]{Nar04} (see also \cite[Corollary 2.3]{LTG24+}).
\end{proof}

The $[a]\in (\Z/m\Z)^\times$ for which $\varepsilon(K,a,m)=1$ form a subgroup of $(\Z/m\Z)^\times$, and for each subgroup $H\leq (\Z/m\Z)^\times$ there are number fields $K$ for which the subgroup 
\begin{equation}\label{eq:sG}
\{[a]\in (\Z/m\Z)^\times : \varepsilon(K,a,m)=1\}
\end{equation}
 equals $H$. As in the introduction, we will denote the subgroup in (\ref{eq:sG}) by $\sG_m(K)$.
We say that $m$ is \textbf{unramified} in $K$ if each rational prime dividing $m$ is unramified in $K$.
When $m$ is unramified in $K$ then $\sG_m(K)=(\Z/m\Z)^\times$. This leads to the following special case of Proposition \ref{prop:K-Dirichlet}.
 
\begin{corollary}
	If $m$ is unramified in $K$, then
	\begin{equation}
		\pi(B;K,a,m)\sim \frac{1}{\deg_K(K(\zeta_m))}\frac{B}{\log(B)}.
	\end{equation}
\end{corollary}

\subsection{Modular curves}

Let $E$ be an elliptic curve over $K$, and let $N$ be a positive integer. Multiplication by $N$ gives a morphism of group schemes
\begin{salign}
m_N: E &\to E\\
P &\mapsto N\cdot P \defeq \underbrace{P+\cdots + P.}_{\text{$N$-times}}
\end{salign}
The $N$-torsion group scheme $E[N]$ is defined as the kernel of $m_N$. Choosing a basis $(P_1,P_2)$ of the group $E[N](\overline{K})$ determines a group isomorphism $E[N](\overline{K})\xrightarrow{\sim} (\Z/N\Z)^2$. This induces an isomorphism $\iota:\Aut(E[N](\overline{K}))\xrightarrow{\sim} \GL_2(\Z/N\Z)$ defined by sending an automorphism $\alpha\in\Aut(E[N](\overline{K}))$ determined by
\begin{salign}
\alpha(P_1)&=a P_1 + c P_2 \\
\alpha(P_2)&= b P_1 + d P_2
\end{salign}
to the matrix 
\begin{equation}
\iota(\alpha)=
\begin{bmatrix}
a & b \\
c & d
\end{bmatrix}.
\end{equation}

Let $S$ be a $K$-scheme, and let $E$ be an elliptic curve over $S$. An $N$-level structure on $E$ is an isomorphism of $S$-group schemes $\iota: (\Z/N\Z)^2_S\xrightarrow{\sim} E[N]$. Call such a pair $(E,\iota)$ an elliptic curve equipped with an $N$-level structure. Two such pairs $(E,\iota)$ and $(E',\iota')$ are equivalent if there is an isomorphism $\varphi:E\xrightarrow{\sim} E'$ such that $\iota=\iota'\circ\varphi$. 

Let $\cY(N)/K$ denote the moduli space of elliptic curves equipped with an $N$-level structure up to equivalence. Let $\cX(N)$ be the smooth compactification of $\cY(N)$. 

Let $G$ be a subgroup of the general linear group $\GL_2(\Z/N\Z)$.
Two $N$-level structures $\iota$ and $\iota'$ are $G$-equivalent if there exists a $g\in G$ such that the following diagram commutes
\begin{center}
\begin{tikzcd}
    & E[N](\overline{K}) \arrow{dl}[swap]{\iota} \arrow{dr}{\iota'} &\\
    (\Z/N\Z)^2 \arrow{rr}{g} & & (\Z/N\Z)^2. \\
\end{tikzcd}
\end{center}
A $G$-level structure on an elliptic curve $E$ is a $G$-equivalence class $[\iota]_G$ of $N$-level structures. Call such a pair $(E,[\iota]_G)$ an elliptic curve equipped with a $G$-level structure. Two such pairs $(E,[\iota]_G)$ and $(E',[\iota']_G)$ are equivalent if there is an isomorphism $\varphi: E\to E'$ such that $\iota$ and $\iota'\circ \varphi$ are $G$-equivalent.

 Let $\cY_G/K$ denote the moduli space of elliptic curves equipped with a $G$-level structure up to equivalence. Let $\cX_G$ be the smooth compactification of $\cY_G$.
 Equivalently, $\cX_G$ can be defined as the quotient stack $[\cX(N)/G]$, where the action of $G$ on $\cX(N)$ (away from cusps) is 
\begin{salign}
G \times \cX(N) &\to \cX(N)\\
(g, (E,\iota)) &\mapsto (E,g\circ \iota).
\end{salign}

An elliptic curve $E$ is said to \emph{admit} a $G$-level structure if it can be equipped with a $G$-level structure.
For example, if $M,N\in \Z_{>0}$ are positive integers and
\[
G(M,MN)\defeq \left\{g\in \GL_2(\Z/MN\Z) : g=
\begin{pmatrix}
\ast & \ast\\
0 & 1
\end{pmatrix}
\text{ and }
g\equiv
\begin{pmatrix}
\ast & 0\\
0 & 1
\end{pmatrix}
\pmod{M}\right\},
\]
then elliptic curves which admit a $G(M,MN)$-level structure are those which have a $K$-rational torsion point $P\in E[MN](K)$ of order $MN$, and a Galois stable cyclic subgroup $C\leq E[M](\overline{K})$ of order $M$ such that $E[M](\overline{K})=\langle NP\rangle \oplus C$.  
 Define $\cX_1(M,MN)\defeq \cX_{G(M,MN)}$ and note that $\cX(N)=\cX_{G(1,N)}$.
 
\begin{remark}
	Note that the definition of $\cX_1(M,MN)$ we use is different from the typical definition of $\cX_1(M,MN)$, as the moduli space of elliptic curves with torsion subgroup $\Z/M\Z \times \Z/MN\Z$. When the base field $K$ contains an $MN$-th root of unity, the set of elliptic curves being parametrized by $\cX_1(M,MN)$ will be the same as those parametrized by the usual definition. The advantage of defining $\cX_1(M,MN)$ in the way we do is that the modular curve is always defined over $\Q$, whereas using the usual definition it is only defined over $\Q(\zeta_{MN})$. A drawback of our definition is that $\cX_1(M,MN)$ will not be geometrically irreducible, but this will not cause us any difficulties. This modular curve is labeled $X_{\text{arith},1}(M,MN)$ in \cite{HL25} and in the L-functions and modular forms database (LMFDB) \cite[\href{https://beta.lmfdb.org/ModularCurve/Q/}{Modular curves over $\Q$}]{lmfdb}.
\end{remark}

Let $\bT$ denote set of $G(M,MN)$ for which $\cX_{G(M,N)}$ has genus zero. It is known that
\begin{salign}
\bT =& \{G(1,N) : 1\leq N\leq 10 \text{ or } N=12\}
 \cup \{G(2,2N) : 1\leq N \leq 4\}\\
&\qquad \cup \{G(3,3), G(3,6), G(4,4), G(5,5)\}.
\end{salign}

For each $G\in \mathbb{T}$ the corresponding modular curve $\cX_G$ is a \textit{weighted projective stack}.

\subsection{Weighted projective stacks}

Let $\A^n$ denote affine $n$-space, and let $\mathbb{G}_m$ denote the multiplicative group scheme. Given a pair of positive integers $(w_0,w_1)$, the \textbf{weighted projective stack} $\cP(w_0,w_1)$ is defined to be the quotient stack
\[
\cP(w_0,w_1)\defeq [(\A^{2}-\{0\})/\mathbb{G}_m],
\]
where $\mathbb{G}_m$ acts on $\A^{n+1}-\{0\}$ by the $(w_0,w_1)$-weighted action
\begin{align*}
\mathbb{G}_m\times(\A^{2}-\{0\}) &\to (\A^{n+1}-\{0\})\\
(\lambda,(x_0,x_1)) &\mapsto  (\lambda^{w_0}x_0, \lambda^{w_1}x_1).
\end{align*}
For example, $\cP(1,1)$ is the usual projective line $\P^1$.

Recall that the elliptic curves $E_1:y^2=x^3+A_1x+B_1$ and $E_2:y^2=x^3+A_2x+B_2$ are isomorphic if and only if there exists a non-zero constant $u\in K^\times$ such that $A_1=u^4 A_2$ and $B_1= u^6 B_2$. It follows  that the compactified moduli space of elliptic curve $\cX(1)=\cX_{G(1,1)}$ is $K$-isomorphic to $\cP(4,6)$, with isomorphism
\begin{salign}
\cX(1) &\xrightarrow{\sim} \cP(4,6)\\
E: y^2=x^3+Ax+B & \mapsto [A:B].
\end{salign}

We now discuss morphisms between weighted projective stacks. Let $\w=(w_0,w_1)$ and $\w'=(w'_0,w'_1)$ be pairs of positive integers. Let $K[x_0,x_1]$ (resp., $K[x'_0,x'_1]$) be a graded $K$-algebra where $x_i$ has weight $w_i$ (reps. $x'_i$ has weight $w'_i$).

 Let $\Psi_{\w',\w}(K)$ denote the set of pairs $(f_1,f_0)\in K[x_0,x_1]$ of weighted homogeneous polynomials satifying the following properties,
 \begin{itemize}
 \item there exists a $\delta\in \Z_{\geq 0}$ such that $\deg(f_i)=\delta\cdot w'_i$, and
 \item the radical ideal $\sqrt{(f_0,f_1)}\subseteq K[x_0,x_1]$ contains the ideal $(x_0,x_1)$. 
 \end{itemize} 
If $f=(f_0,f_1)\in \Psi_{\w',\w}(K)$ then the ring homomorphism
\begin{align*}
K[x'_0,x'_1]&\to K[x_0,x_n]\\
x_i &\mapsto f_i
\end{align*}
induces a morphism of weighted projective stacks
\[
\psi_f:\cP(w_0,w_1)\to \cP(w'_0,w'_1)
\]
(see \cite[Lemma 4.1]{BN22}).
We call $\delta$ the \textbf{weighted degree} of the morphism $\psi_f$.
 There is a natural weighted action of $\G_m(K)$ on $\Psi_{\w',\w}(K)$ defined as
\begin{align*}
\G_m(K)\times \Psi_{\w',\w}(K)&\to \Psi_{\w',\w}(K)\\
(\lambda, (f_i)) & \mapsto (\lambda^{w_0} f_0, \lambda^{w_1} f_1).
\end{align*} 
Let $\Hom_K(\cP(\w),\cP(\w'))$ denote the set of morphisms from $\cP(\w)$ to $\cP(\w)'$ defined over $K$.

\begin{proposition}[\protect{\cite[Lemma 4.1]{BN22}}]\label{prop:WProjMorphisms}
Using the notation above, there is a bijection
\begin{align*}
\Phi_{\w',\w}(K)/\G_m(K) &\to \Hom_K(\cP(\w),\cP(\w'))\\
f=(f_0,f_1) &\mapsto \psi_f.
\end{align*}
\end{proposition}

\subsection{Heights}\label{subsection:heights}
We define a height function on elliptic curves over a number field by defining a height on the set of $K$-points $\cP(4,6)(K)$.

For $x=(x_0,x_1)\in K^{2}-\{0\}$, let $|x_i|_{(w_0,w_1),v}\defeq |\pi_v|_v^{\lfloor v(x_i)/w_i\rfloor}$ and set
\begin{equation}
 |x|_{{(w_0,w_1)},v}\defeq \begin{cases}
 \max_i\{|x_i|_{(w_0,w_1),v}\} & \text{ if } v\in \Val_0(K),\\
 \max_i\{|x_i|_v^{1/w_i}\} & \text{ if } v\in \Val_\infty(K).
 \end{cases}
 \end{equation}
 
  The \textbf{height} of a point $x=[x_0:\cdots:x_n]\in \cP(\w)(K)$ is defined as
\begin{equation}
\Ht_{(w_0,w_1)}(x)\defeq \prod_{v\in \Val(K)} |(x_0,\dots,x_n)|_{(w_0,w_1),v}.
\end{equation}

The \textbf{naive height} of an elliptic curve $E$ over the number field $K$ is
\begin{equation}\label{eq:naive_height}
\Ht(E)\defeq \Ht_{(4,6)}([a:b])^{12}.
\end{equation}

\subsection{Isogenies}\label{subsec:isogenies}

An \textbf{isogeny} $\phi$ between elliptic curves $E_1$ and $E_2$ is a finite morphism of group scheme $\phi: E_1\to E_2$.
Two elliptic curves $E_1$ and $E_2$ are \textbf{isogenous} over $K$ if there exists an isogeny defined over $K$ between them. Let $E[\phi]$ denote the kernel of $\phi$. The \textbf{degree} of an isogeny is its degree as a morphism of algebraic curves. As we are working over number fields, the degree coincides with the cardinality of the $\overline{K}$-points of the kernel $E[\phi](\overline{K})$. If $\phi:E_1\to E_2$ is an isogeny of degree $N$, then there exists a unique isogeny $\widehat{\phi}:E_2\to E_1$ with the property that $\widehat{\phi}\circ\phi=m_N$ \cite[Theorem 6.1]{Sil09}; this isogeny is called the \textbf{dual isogeny}.

\begin{example}
The multiplication-by-$N$ morphism $m_N$ is an isogeny of degree $N^2$.
\end{example}

\begin{proposition}[\protect{\cite[Proposition 4.12]{Sil09}}]\label{proposition:subgroup_isogeny}
If $E$ is an elliptic curve and $\Phi$ is a finite subgroup of $E$, then there exists a unique elliptic curve $E'$ and isogeny $\phi:E\to E'$ with kernel $E[\phi]=\Phi$.
\end{proposition}


Let $G=G(M,MN)\in \mathbb{T}$, let $\Phi$ be a subgroup of $\Z/M\Z\times \Z/MN\Z$, and let $d=\# \Phi$. By Proposition \ref{proposition:subgroup_isogeny}, for each elliptic curve $E$ which admits a $G$-level structure there is an isogeny $\phi:E\to E'$ with kernel $\Phi$. This induces a morphism of stacks
\begin{salign}\label{eq:isogeny_morphism}
\phi_\Phi:\cX_G &\to \cX(1)\\
(E,[\iota]_G) &\mapsto E'.
\end{salign}
The modular curve $\cX_{G}$ is isomorphic to a weighted projective stack $\cP(w_0,w_1)$ with weights $(w_0,w_1)$ given in Table \ref{tab:elliptic_count}.
   Therefore the morphism (\ref{eq:isogeny_morphism}) is equivalent to a morphism of weighted projective stacks,
 \begin{salign}
    \phi_\Phi:\cP(w_0,w_1) &\to \cP(4,6)\\
    (A,B) &\mapsto (f_{G,4}(A,B),f_{G,6}(A,B)),
\end{salign}
where, by Proposition \ref{prop:WProjMorphisms}, $f_{G,4}$ and $f_{G,6}$ are weighted homogeneous polynomials with degrees $\deg(f_{G,4})=4n$ and $\deg(f_{G,6})=6n$ for some  $n\in \Z_{>0}$. 
  Let $\Delta_{G,\Phi}(A,B)$ denote the discriminant the elliptic curve  
  \begin{equation}
  E: y^2=x^3+f_{G,4}(A,B)x+f_{G,6}(A,B).
  \end{equation}

When $\Phi=\{O\}$, we have that $\phi_{O}: \cX_G\to \cX(1)$ is the morphism that forgets the level structure. In these cases $f_{G,4}(A,B)$ and $f_{G,6}(A,B)$ can be found in Table \ref{tab:Weierstrass_families}. In most cases these weighted homogeneous polynomials can be found by starting with the universal family of elliptic curves with $G$-level structure (see, e.g., \cite[Table 3]{Kub76}), and then transforming it into short Weierstrass form and finally homogenizing with respect to the proper weights. However, in the cases of $G(3,3)$, $G(3,6)$, $G(4,4)$, and $G(5,5)$ such universal families do not appear to be in the literature. In these cases one can find a model for the universal family by finding a suitable twist of an elliptic curve with the correct $j$-invariant (the relevant $j$-invariants can be found on the LMFDB \cite{lmfdb}).

For general $\Phi$, defining equations can be computed from those in the $\Phi=\{O\}$ case by using Velu's formula. Velu's formula \cite{Vel71} gives an algorithm for computing the isogeny in Proposition \ref{proposition:subgroup_isogeny}. 

Table \ref{tab:Discriminants} lists the discriminants $\Delta_{G,\Phi}(A,B)$ for each $G\in \bT$ when $\#\Phi$ is one or prime.



 \begin{example}
 	Every elliptic curve with a torsion point $P$ of order $5$ can be written in Tate normal form as
 	\begin{equation}\label{eq:Tate_normal_5}
 		y^2 + (-A + B)xy - AB^2y = x^3 - ABx^2,
 	\end{equation}
 	with $(A,B)\in \cP(1,1)(K)$, and where the $5$-torsion point is $P=(0,0)$. Rewritting (\ref{eq:Tate_normal_5}) in short Weierstrass form allows us to compute explicit equations for the morphism that forgets the level structure
 	\begin{salign}
 		\phi_{O}: \cP(1,1) &\to \cP(4,6)\\
 		(A,B) &\mapsto (g_{G(1,5), 4}, g_{G(1,5),6}),
 	\end{salign}
 	where
 	 \begin{salign}
 		 g_{G(1,5), 4} &=-27A^4 + 324A^3B - 378A^2B^2 - 324AB^3 - 27B^4 \\
 		g_{G(1,5),6} &=54A^6 - 972A^5B + 4050A^4B^2 + 4050A^2B^4 + 972AB^5 + 54B^6.
 	\end{salign}
 	The discriminant is
 	\begin{equation}
 		\Delta_{G(1,5),O}(A,B)=A^5B^5(A^2 - 11AB - B^2).
 	\end{equation}

 	Using Velu's formula \cite{Vel71}, we compute
 	\begin{salign}
 		\phi_{\langle P\rangle}: \cP(1,1) &\to \cP(4,6)\\
 		(A,B) &\mapsto (f_{G(1,5), 4}, f_{G(1,5),6}),
 	\end{salign}
 	where
 	\begin{salign}
 		 f_{G(1,5), 4} &=-27A^4 - 6156A^3B - 13338A^2B^2 + 6156AB^3 - 27B^4 \\
 		f_{G(1,5),6} &=54A^6 - 28188A^5B - 540270A^4B^2 - 540270A^2B^4 + 28188AB^5 + 54B^6.
 	\end{salign}
 	The discriminant is
 	\begin{equation}
 		\Delta_{G(1,5),\langle P\rangle}(A,B)=AB(A^2 - 11AB - B^2)^5.
 	\end{equation}
 \end{example}

 \begin{graytext}


\end{graytext}

\subsection{Selmer groups}



Let $\phi:E\to E'$ be an isogeny of elliptic curves over $K$.

There is a short exact sequence of $\cG_K$-modules
\begin{equation}
\begin{tikzcd}
0\rar & E[\phi] \rar & E \ar[r,"\phi"] & E' \rar  & 0.
\end{tikzcd}
\end{equation}
Taking Galois cohomology we form the long exact sequence
\begin{equation}\label{eq:long_exact}
\begin{tikzcd}
  0\rar & E[\phi](K) \rar & E(K) \ar[r,"\phi"]
             \ar[draw=none]{d}[name=X, anchor=center]{}
    & E'(K) \ar[rounded corners,
            to path={ -- ([xshift=2ex]\tikztostart.east)
                      |- (X.center) \tikztonodes
                      -| ([xshift=-2ex]\tikztotarget.west)
                      -- (\tikztotarget)}]{dll}[at end, swap]{\delta_\phi}\\      
 & H^1(K,E[\phi]) \rar & H^1(K,E) \ar[r,"\phi"] & H^1(K,E'),
\end{tikzcd}
\end{equation}
where $\delta_\phi$ is the connecting homomorphism.

 From (\ref{eq:long_exact}) we obtain the short exact sequence
\begin{equation}
\begin{tikzcd}
0\rar & E'(K)/\phi(E(K)) \ar[r,"\delta_\phi"] & H^1(K,E[\phi]) \ar[r] & H^1(K,E)[\phi] \rar  & 0.
\end{tikzcd}
\end{equation}
Similarly there is an exact sequence
\begin{equation}
\begin{tikzcd}
0\rar & E'(K_v)/\phi(E(K_v)) \ar[r,"\delta_{\phi,v}"] & H^1({K_v},E[\phi]) \ar[r] & H^1({K_v},E)[\phi] \rar  & 0.
\end{tikzcd}
\end{equation}

The inclusions $\cG_{K_v}\hookrightarrow \cG_K$ and $E(K^{\sep})\hookrightarrow E(K_v^{\sep})$ induce restriction maps, leading to the following commutative diagram:
\begin{equation*}
\begin{tikzcd}
0\rar & E'(K)/\phi(E(K)) \ar[r,"\delta_\phi"] \dar & H^1(K,E[\phi]) \ar[r]\dar 
& H^1(K,E)[\phi] \rar  \dar & 0\\
0\rar & \prod\limits_{v}E'(K_v)/\phi(E(K_v)) \ar[r,"\delta_{\phi,v}"] & \prod\limits_{v}H^1({K_v},E[\phi]) \ar[r] & \prod\limits_{v} H^1({K_v},E)[\phi] \rar  & 0,
\end{tikzcd}
\end{equation*}
where the products are over all place $v\in \Val(K)$.

Define subgroups
\begin{equation}
H^1_\phi({K_v},E[\phi])\defeq \delta_{\phi,v}\left(E'(K_v)/\phi(E(K_v))\right)\subset H^1({K_v}, E[\phi]).
\end{equation}
Let 
\begin{equation}
	\res_v: H^1(K,E[\phi])\to H^1(K_v, E[\phi])
\end{equation}
denote the natural restriction maps.
The \textbf{$\phi$-Selmer group} of $E$ is
\begin{equation}
	\Sel_\phi(E)\defeq \{x\in H^1(K,E[\phi]) : \res_v(x)\in H^1_\phi({K_v},E[\phi]) \text{ for all } v\in \Val_0(K)\}.
\end{equation}

The \textbf{$N$-Selmer group}, denoted $\Sel_N(E)$, is the $\phi$-Selmer group when $\phi=m_N$ is the multiplication-by-$N$ isogeny.

Let $E$, $E'$, and $E''$ be elliptic curves and let $\phi,\phi_1,\phi_2$ be isogenies such that the following diagram commutes:
\begin{equation}
\begin{tikzcd}[row sep=4em, column sep=2em]
    E \arrow[rr, "\phi"] \arrow[swap, dr, "\phi_1"] & & E''  \\
    & E' \arrow[swap, ur, "\phi_2"]
\end{tikzcd}
\end{equation}

There is a short exact sequence of $\cG_K$-modules
\begin{equation}
\begin{tikzcd}
0\rar & E[\phi_1] \rar & E[\phi] \ar[r,"\phi_1"] & E'[\phi_2] \rar  & 0.
\end{tikzcd}
\end{equation}
Taking Galois cohomology and using the definition of the Selmer groups, we arrive at the long exact sequence
\begin{equation}\label{eq:long_exact_phi1-phi-phi2}
\begin{tikzcd}
  0\rar & E[\phi_1](K) \rar & E[\phi](K) \ar[r,"\phi_1"]
             \ar[draw=none]{d}[name=X, anchor=center]{}
    & E'[\phi_2](K) \ar[rounded corners,
            to path={ -- ([xshift=2ex]\tikztostart.east)
                      |- (X.center) \tikztonodes
                      -| ([xshift=-2ex]\tikztotarget.west)
                      -- (\tikztotarget)}]{dll}[at end, swap]{\delta_{\phi_1}}\\      
 & \Sel_{\phi_1}(E) \rar & \Sel_\phi(E) \ar[r,"\phi"] & \Sel_{\phi_2}(E').
\end{tikzcd}
\end{equation} From this we obtain the exact sequence
\begin{equation}\label{eq:medium_exact_phi1-phi-phi2}
0
\rightarrow \frac{E'[\phi_2](K)}{\phi(E[\phi](K))} 
\rightarrow \Sel_{\phi_1}(E)
\rightarrow \Sel_\phi(E)
\xrightarrow{\phi_1} \Sel_{\phi_2}(E').
\end{equation}

\begin{proposition}\label{prop:Selmer_size_comparison}
Let $\phi_1:E\to E'$ and $\phi_2:E'\to E''$ be isogenies of elliptic curves over a number field $K$, and let $\phi=\phi_2\circ\phi_1$. Then
\begin{equation}
	\#\Sel_{\phi_1}(E)\leq \deg(\phi_2)\cdot \# \Sel_{\phi}(E).
	\end{equation}
\end{proposition}

\begin{proof}
	Since (\ref{eq:medium_exact_phi1-phi-phi2}) is an exact sequence of finite groups,
	\begin{equation}
		\#\Sel_{\phi_1}(E)\leq \#\left(\frac{E'[\phi_2](K)}{\phi_1(E[\phi](K))} \right) \cdot \# \Sel_{\phi}(E)\leq \deg(\phi_2)\cdot \#\Sel_{\phi}(E).
	\end{equation}
\end{proof} 

If $\phi_1$ and $\phi_2$ are dual isogenies of degree $N$, then $\phi$ is the multiplication-by-$N$ isogeny $m_N$. By Proposition \ref{prop:Selmer_size_comparison},
\begin{equation}\label{eq:N-Selmer_inequality}
	\# \Sel_N(E) \geq \frac{\#\Sel_{\phi_1}(E)}{N} \geq \frac{\#\Sel_{\phi_1}(E)}{N \# \Sel_{\widehat{\phi}_1}(E)}.
\end{equation}
The reason for dividing by $\# \Sel_{\widehat{\phi}_1}(E)$ is that the ratio $\#\Sel_{\phi_1}(E)/\# \Sel_{\widehat{\phi}_1}(E)$ has nice properties, which we make key use of in this article. Determining the average behavior of this ratio will give us a lower bound for the average size of the $N$-Selmer group.

\section{Selmer ratios and Tamagawa numbers}\label{sec:Selmer_ratios}

In this section, we study ratios of certain Selmer groups and see how they are related to ratios of certain Tamagawa numbers.

Define the \textbf{Selmer ratio} of $E$ to be
\begin{equation}\label{eq:Selmer_ratio}
\mathcal{S}(E/E')\defeq \frac{\# \Sel_\phi(E)}{\# \Sel_{\widehat{\phi}}(E')}.
\end{equation}

\begin{proposition}
The Selmer ratio decomposes into a product of local factors,
\begin{equation}\label{eq:local_Selmer_ratio}
\cS_v(E/E')\defeq \frac{\# H^1_\phi(K_v,E[\phi])}{\deg(\phi)},
\end{equation}
as follows:
\begin{equation}
\cS(E/E')=\prod_{v\in \Val(K)} \cS_v(E/E').
\end{equation}
\end{proposition} 

\begin{proof}
The groups $H^1_\phi(K_v,E[\phi])$ and $H^1_{\widehat{\phi}}(K_v,E'[\widehat{\phi}])$ are exact annihilators with respect to the Tate local pairing between $H^1(K_v,E[\phi])$ and $H^1(K_v,E'[\widehat{\phi}])$ \cite[Lemma 2.4]{Fis03}. Therefore the desired product formula is a special case of \cite[Theorem 2.19]{DDT94}.
\end{proof}

\begin{proposition}\label{prop:local_selmer_ratio_bound}
We have the following bounds for the local factors:
\begin{equation}
	\cS_v(E/E')\leq 
	\begin{cases}
		\deg(\phi)^2 & \text{ if } v\in \Val_0(K),\\
		2/\deg(\phi) & \text{ if } v\in \Val_\infty(K).
	\end{cases}
\end{equation}
\end{proposition}

\begin{proof}
From the definitions of $\cS_v(E/E')$ and $H_\phi^1(K_v, E[\phi])$, it suffices to show that
\begin{equation}
\#\left(\frac{E'(K_v)}{\phi(E(K_v))}\right)\leq 			 	\begin{cases}
		\deg(\phi)^3 & \text{ if } v\in \Val_0(K),\\
		2 & \text{ if } v\in \Val_\infty(K).
	\end{cases}
\end{equation}

 The case $v\in \Val_\infty(K)$ is part of \cite[Lemma 3.11]{Sch96}.
 
 Now consider the case $v\in \Val_0(K)$.
 As the image of the multiplication-by-$\deg(\phi)$ map $m_{\deg(\phi)}(E'(K_v))$ is contained in $\phi(E(K_v))$,
\begin{equation}
\#\left(\frac{E'(K_v)}{\phi(E(K_v))}\right)\leq\#\left(\frac{E'(K_v)}{m_{\deg(\phi)}(E(K_v))}\right).
\end{equation}
Let $\p_v$ denote the prime ideal corresponding to $v$. By \cite[Proposition 3.9]{Sch96},
\begin{equation}
	\#\left(\frac{E'(K_v)}{m_{\deg(\phi)}(E(K_v))}\right)=N_{K/\Q}(\p_v)^{v(\deg(\phi))}\#E[\deg(\phi)](K_v)\leq \deg(\phi)^3.
\end{equation}
\end{proof}

\begin{question}
	If $E[\phi]$ is a constant group scheme, then it can be shown that, for any place $v$,
    \begin{equation}
        \# H^1(K_v,E[\phi])=O_{\deg(\phi)}(1).
    \end{equation} 
	Is this true for any isogeny $\phi$, without the restriction that $E[\phi]$ is a constant group scheme?
\end{question}



Let $E_0(K_v)$ denote the identity component of $E(K_v)$. 
The \textbf{local Tamagawa number} of $E$ at $v$ is the index 
\begin{equation}
c_v(E)\defeq [E(K_v):E_0(K_v)].
\end{equation}

\begin{proposition}\label{prop:tamagawa_ratio}
If $v(\deg(\phi)) = 0$, then
\begin{equation}
\cS_v(E/E')=\frac{c_v(E')}{c_v(E)}.
\end{equation}
\end{proposition}

\begin{proof}
This is a special case of \cite[Lemma 3.8]{Sch96}.\footnote{See the beginning of the second paragraph after the proof of \cite[Lemma 3.8]{Sch96}.}
\end{proof}

 Let $\Delta(E)$ denote the minimal discriminant of $E$ over $K$.
 Using Tate's algorithm \cite{Tat75}, one can prove the following facts about local Tamagawa numbers. 

\begin{proposition}\label{prop:tamagawa_numbers}
Let $v\in \Val_0(K)$ be a finite place.
\begin{enumerate}[label=(\alph*)]
\item If $E$ has split-multiplicative reduction at $v$ then $c_v(E)=v(\Delta(E))$.
\item If $E$ has non-split multiplicative reduction at $v$ then $c_v\leq 2$.
\item If $E$ has additive reduction at $v$ then $c_v\leq 4$.
\item If $E$ has good reduction at $v$, then $c_v=1$.
\end{enumerate}
\end{proposition}

By Proposition \ref{prop:tamagawa_ratio} and Proposition \ref{prop:tamagawa_numbers}, if $E$ has good reduction at a finite place $v$, then $\cS_v(E/E')=1$. In particular, we have that the Selmer ratio is a finite product
\begin{equation}
\cS(E/E')=\prod_{v(\deg(\phi)\cdot \Delta(E))>0} \cS_v(E/E') \prod_{v\in \Val_\infty(K)} \cS_v(E/E').
\end{equation} 

Every isogeny can be decomposed into a composition of isogenies of prime degree. We have the following result concerning Tamagawa ratios $c_v(E')/c_v(E)$ when $E$ has multiplicative reduction at $v$.

\begin{proposition}\label{prop:tamagawa_ratio_values}
Let $\phi:E\to E'$ be an isogeny of prime degree $p$, and suppose that $E$ has multiplicative reduction at a place $v\in \Val_0(K)$ not above $p$. If $p\neq 2$, then
\begin{equation}
\frac{c_v(E')}{c_v(E)}=
\begin{cases}
\frac{v(\Delta(E'))}{v(\Delta(E))} & \text{ if $E$ has split multiplicative reduction at $v$,}\\
1 & \text{ if $E$ has non-split multiplicative reduction at $v$}.
\end{cases}
\end{equation}
If $p=2$, then
\begin{equation}
\frac{c_v(E')}{c_v(E)}=
\begin{cases}
\frac{v(\Delta(E'))}{v(\Delta(E))} & \text{ if $ 2\nmid\gcd(v(\Delta(E)),v(\Delta(E')))$ or $E$ has split mult.\@ reduction at $v$},\\
1 & \text{ if $2|\gcd(v(\Delta(E)),v(\Delta(E')))$ and $E$ has non-split mult.\@ reduction at $v$}.
\end{cases}
\end{equation}
\end{proposition}

\begin{proof}
This is part of \cite[Table 1]{DD15}.
\end{proof}

Let $\phi$ and $v$ be as in Proposition \ref{prop:tamagawa_ratio_values}, and let $n\in\{-1,0,1\}$. We say that $E$ has \textbf{$n$-multiplicative reduction} at $v$ with respect to $\phi$ if $E$ has multiplicative reduction at $v$ and
$n=\log_p\left(c_v(E')/c_v(E)\right)$. 

 \section{Counting elliptic curves with prescribed torsion and local conditions}\label{sec:local_conditions}

In this section, we explain how results of \cite{Phi24+} can be extended to count elliptic curves over number fields with a prescribed torsion subgroup and prescribed local conditions.

Let $G\leq \GL_2(\Z/N\Z)$ be such that the modular curve $\cX_G$ is isomorphic to the weighted projective stack $\cP(w_0,w_1)$.
 For a finite place $v\in \Val_0(K)$, and a subset $\Omega\subseteq \cP(w_0,w_1)(K_v)$, let 
 \begin{equation}
 \Omega^{\aff}\defeq \{(x_0,\dots,x_n)\in K_v^{n+1} :  [x_0:\dots:x_n]\in \Omega\}
 \end{equation}
 denote the subset of the affine cone of $\cP(w_0,w_1)$ above $\Omega$, and let $\partial\Omega^{\aff}$ denote the boundary of $\Omega^{\aff}$. 
We will say that $E$ is contained in $\Omega_v$, and write $E\in \Omega_v$, if the point $(a,b)\in \cP(w_0,w_1)(K_v)$ corresponding to $E$ is contained in $\Omega_v$.

Let $\cE_G(K)$ denote the set of elliptic curves defined over $K$ which admit a $G$-level structure, and for each $B\in \R_{\geq 0}$ define the set of elliptic curves in $\cE_G(K)$ of bounded height,
\begin{equation}
\cE_G(K;B)\defeq \{E\in \cE_G(K) : \Ht(E)\leq B\}. 
\end{equation}
Define the counting function
\begin{equation}
\cN_G(K,\Omega_v;B)\defeq \#\{E\in \cE_G(B;K) : E\in \Omega_v \}.
\end{equation}

\begin{proposition}\label{prop:elliptic_count_general}
Let $G\leq \GL_2(\Z/N\Z)$, and suppose that the modular curve $\cX_G$ is isomorphic to the weighted projective stack $\cP(w_0,w_1)$ over the number field $K$. 
Let $\delta(G)$ denote the weighted degree of the morphism $\phi_{O}:\cX_G\to\cX(1)$ that forgets the level structure, and suppose that $\gcd(w_0w_1, \delta(G))=1$. Define
\begin{equation}
 \alpha(G)\defeq\frac{w_0+w_1}{12 \delta(G)}\qquad \text{ and }\qquad
 \beta(G)\defeq \frac{\deg_{\Q}(K)-\min\{w_0,w_1\}}{12\delta(G)\deg_{\Q}(K)}.
 \end{equation}
Let $v\in \Val_0(K)$ be a finite place not above $2$ or $3$, and let $\Omega_v\subseteq\cP(w_0,w_1)(K_v)$ be a subset for which $\Omega_v^{\aff}$ is measurable and whose boundary has measure $m_v(\partial \Omega_v^{\aff})=0$. 
Then there exists an explicit constant $\kappa_G$, depending only on $G$ and $K$, such that
\begin{equation}
\cN_G(K,\Omega_v; B) 
 = m_v\left(\Omega_v^{\aff}\cap \cO_{K,v}^{2}\right)\kappa_G B^{\alpha(G)}
  + O\left(B^{\beta(G)}\right).
\end{equation}
\end{proposition}

\begin{proof}
	This is a slight generalization of \cite[Theorem 1.2.1]{Phi24+}, which proved this result without the local condition $\Omega_v$. The same proof works in this case, noting that \cite[Theorem 1.1.1]{Phi24+} counts points with (finitely many) local conditions.
\end{proof}

We use Proposition \ref{prop:elliptic_count_general}  to count elliptic curves with a prescribed torsion subgroup and a specified local condition.
For $v\in \Val_0(K)$ a finite place, let
\begin{equation}
\star\in\{\text{good},\ \text{additive},\ \text{multiplicative},\ \text{split},\ n\text{-multiplicative}\}.
\end{equation} 
Define the counting function
\begin{equation}
\cN_G(K,v,\star; B)\defeq \#\{E\in \cE_G(K;B) : E \text{ has $\star$ reduction at $\p_v$}\}.
\end{equation}

\begin{proposition}\label{prop:elliptic_count}
 Maintain the notation and assumptions of Proposition \ref{prop:elliptic_count_general}. Then there exists explicit constants $\kappa_\star(G,v)$ and $\kappa_G$ such that 
 \begin{equation}
\cN_G(K,v,\star; B)
=\left(\frac{\kappa_{\star}(G,v)  \kappa_G}{1-N_{K/\Q}(\p_v)^{-10}}\right) B^{\alpha(T)} 
+ O\left(B^{\beta(T)}\right).
\end{equation} 
 For $G\in \bT$ the values of $\alpha(G)$ and $\beta(G)$ are given in Table \ref{tab:elliptic_count}, and the values of $\kappa_\star(G,v)$ are given in Table \ref{tab:LocalConditions} and Table \ref{tab:isogeny_local_densities}. An expression for the constant $\kappa_G$ is given in \cite[Section 5, Theorem 1.2.1]{Phi24+}. 
 \end{proposition}


Proposition \ref{prop:elliptic_count} is an application of Proposition \ref{prop:elliptic_count_general}. The computations for the local densities in Table \ref{tab:LocalConditions} and Table \ref{tab:isogeny_local_densities} are similar to the examples given in the the proof of \cite[Theorem 1.1.2]{Phi25}, but less straightforward. To illustrate what is involved in these computations, we give a couple examples.

\begin{example}[Counting elliptic curves with $3$-torsion and split multiplicative reduction]

Every elliptic curve with a point of order $3$ has a model of the form
\begin{equation}\label{eq:3-torsion_model}
E: y^2 + Axy + By = x^3,
\end{equation}
where the point of order $3$ is $(0,0)$. This can be transformed into the short Weierstrass model
\begin{equation}\label{eq:3-torsion_short_weierstrass}
E: y^2 = x^3 - 27(A^4 + 216AB) X + 54(A^6 - 540A^3B - 5832B^2).
\end{equation}
As $\cX_1(3)\cong \cP(1,3)$ and $\cX(1)\cong \cP(4,6)$, this defines a morphism 
	\begin{align*}
		\phi_{O} : \cP(1,3) &\to \cP(4,6)\\
				[A:B] &\mapsto [-27(A^4 + 216AB): 54(A^6 - 540A^3B - 5832B^2)],
	\end{align*}
	which is the morphism that forgets the level structure. 
	
	By Tate's algorithm, the elliptic curve (\ref{eq:3-torsion_model}) has multiplicative reduction at a finite place $v\in \Val_0(K)$ if and only if exactly one of the following is true:
	\begin{itemize}
		\item $v(B)\geq 1$,
		\item $v(A^3-27B)\geq 1$.
	\end{itemize}
	In the first case, $E$ will always have split multiplicative reduction. In the second case, $E$ will have split multiplicative reduction if and only if $-3$ is a square in the residue field $\F_{q_v}$ at $v$.
	
We now count the number of pairs $(a,b)\in \F_{q_v}\times \F_{q_v}$ with the following properties:
\begin{itemize}
	\item There exists $(A,B)\in \F_{q_v}\times \F_{q_v}$ with 
	\begin{equation}
	(a,b)=(-27(A^4 + 216AB), 54(A^6 - 540A^3B - 5832B^2)),
	\end{equation}
	 and
	\item exactly one of the following is true:
	\begin{itemize}
		\item[$\circ$] $v(B)\geq 1$,
		\item[$\circ$] we have $v(A^3-27B)\geq 1$ and $-3$ is a quadratic residue in $\F_{q_v}$.
	\end{itemize}
\end{itemize} 
The number of such pairs can be shown to be
\begin{equation}
	\begin{cases}
		2q_v-2 & \text{ if } q_v\equiv 1\pmod{3},\\
		q_v-1 & \text{ if } q_v \equiv 2 \pmod{3}.
	\end{cases}
\end{equation}
For each such pair $(a,b)$ consider a lift to $\cO_{K,v}^2$, which we also denote by $(a,b)$. Apply Proposition \ref{prop:elliptic_count_general} with the local condition
\begin{equation}
	\Omega_v(a,b) = \{ [x_0:x_1] \in \cP(w_0,w_1)(K_v) : [x_0:x_1]=[a:b]\}.
\end{equation}
To compute the measure
\begin{equation}
	m_v\left(\Omega_{v}^{\aff}(a,b)\cap \cO_{K,v}^2 \right),
\end{equation}
define the sets
\begin{equation}
\Omega_{v,t}^{\aff}(a,b)
 = \left\{(x_0,x_1)\in \cO_{K,v}^{2} : |x_0-q^{4t}a|_v\leq \frac{1}{q^{4t+1}},\ |x_1-q^{6t}b|_v\leq \frac{1}{q^{6t+1}}\right\},
\end{equation}
which give a partition
\begin{equation}
	\Omega_{v}^{\aff}(a,b)\cap \cO_{K,v}^2 
	= \bigsqcup_{t\in \Z_{\geq 0}} \Omega_{v,t}^{\aff}(a,b).
\end{equation}
Since $\Omega_{v,t}^{\aff}(a,b)$ has $v$-adic measure
\begin{equation}
	m_v\left(\Omega_{v,t}^{\aff}(a,b)\right)=\frac{1}{q_v^{10t+2}},
\end{equation}
the $v$-adic measure of $\Omega_{v}^{\aff}(a,b)\cap \cO_{K,v}^2$ is 
\begin{equation}
	m_v\left(\Omega_{v}^{\aff}(a,b)\cap \cO_{K,v}^2 \right) = \sum_{t\geq 0} \frac{1}{q^{10t+2}}=\frac{1}{q_v^2} \frac{q_v^{10}}{q_v^{10}-1}.
\end{equation} 

 Summing over all possible pairs $(a,b)$ shows that the number of elliptic curves over $K$, with a point of order $3$, split multiplicative reduction at $v$, and of height bounded by $B$ is 
\begin{align}
\frac{q_v^{10}}{q_v^{10}-1} \kappa_{split}(G(1,3),v) \kappa_G B^{1/3}+O\left(B^{\frac{1}{3}-\frac{1}{12\deg_{\Q}(K)}}\right),
\end{align}
where 
\[
\kappa_{split}(G(1,3),v) = 
\begin{cases}
\frac{2q_v-2}{q_v^2}  \text{ if } q_v\equiv 1 \pmod{3},\\
\frac{q_v-1}{q_v^2}  \text{ if } q_v\equiv 2 \pmod{3}.
\end{cases}
\]
\end{example}

\begin{example}[Counting elliptic curves with $4$-torsion and $1$-reduction at the isogeny $\phi$]

Every elliptic curve with a point of order $4$ has a model of the form
\begin{equation}\label{eq:4-torsion_model}
E: y^2 + Bxy -ABy = x^3 - Ax^2,
\end{equation}
where $P=(0,0)$ is a point of order four. Let $\Phi=\langle 2P\rangle$ be the order two subgroup of $\langle P\rangle\cong \Z/4\Z$. By Proposition \ref{proposition:subgroup_isogeny}, there exists an elliptic curve $E'$ and a degree two isogeny
\begin{equation}
	\phi: E \to E'
\end{equation}
with kernel $\Phi$. Using Velu's formula, a model for $E'$ is
\begin{equation}\label{eq:4-torsion_velu_model}
E': y^2 + Bxy -AB^2y = x^3 - Ax^2 - 5A^2x + (-3A^3 - A^2B^2).
\end{equation}
These elliptic curves have discriminant
\begin{equation}
	\Delta(E) = A^4B^2(16A+B^2) \quad \text{ and }\quad \Delta(E') = A^2B^4(16A+B^2)^2.
\end{equation}
The models (\ref{eq:4-torsion_model}) and (\ref{eq:4-torsion_velu_model}) can be transformed into the short Weierstrass models
\begin{equation*}\label{eq:4-torsion_sw}
E: y^2  = x^3 - 27(16A^2 + 16AB^2 + B^4) x - 54(64A^3 - 120A^2B^2 - 24AB^4 - B^6)
\end{equation*}
and
\begin{equation*}\label{eq:4-torsion_velu_sw}
E': y^2  = x^3 - 27(256A^2 + 16AB^2 + B^4) x + -54(4096A^3 + 384A^2B^2 - 24AB^4 - B^6).
\end{equation*}

We have $\cX_1(4)\cong \cP(1,2)$ and $\cX(1)\cong \cP(4,6)$. As in Subsection \ref{subsec:isogenies}, the given models for $E$ and $E'$ define morphisms
\begin{align*}
		\phi_{O} : \cP(1,2) &\to \cP(4,6)\\
				[A:B] &\mapsto [-27(16A^2 + 16AB^2 + B^4): -54(64A^3 - 120A^2B^2 - 24AB^4 - B^6)],
\end{align*}
and
\begin{align*}
		\phi_{\langle 2P\rangle} : \cP(1,2) &\to \cP(4,6)\\
				[A:B] &\mapsto [-27(256A^2 + 16AB^2 + B^4) : -54(4096A^3 + 384A^2B^2 - 24AB^4 - B^6)],
\end{align*}
respectively.

	By Tate's algorithm and Proposition \ref{prop:tamagawa_ratio_values}, the elliptic curve (\ref{eq:4-torsion_model}) has $1$-multiplicative reduction at a finite place $v\in \Val_0(K)$ (\emph{not} above $2$ or $3$) with respect to $\phi$ if and only if one of the following is true:
	\begin{itemize}
		\item $v(B)>0$, $v(A)=0$, and $-A$ is a square in $\F_{q_v}$,
		\item $v(16A+B^2)>0$ and $v(A)=0$. 
	\end{itemize}
	
	
Now we count the number of pairs $(a,b)\in \F_{q_v}\times \F_{q_v}$ with the following properties:
\begin{itemize}
	\item There exists $(A,B)\in \F_{q_v}\times \F_{q_v}$ with 
	\begin{equation*}
	(a,b)=(-27(16A^2B^2 + 16AB^3 + B^4), -54(64A^3B^3 - 120A^2B^4 - 24AB^5 - B^6)),
	\end{equation*}
 and one of the following is true:
 \begin{itemize}
 	\item[$\circ$] $v(B)>0$, $v(A)=0$, and $-A$ is a square in $\F_{q_v}$,
 	\item[$\circ$] $v(16A+B^2)>0$ and $v(A)=0$.
 \end{itemize}
\end{itemize} 
A computation shows that the number of such pairs is $3(q_v-1)/2$.

For each such pair $(a,b)$ choose a lift to $\cO_{K,v}^2$, which we also denote by $(a,b)$. 
As in the previous example, apply Proposition \ref{prop:elliptic_count_general} with the local condition
\begin{equation}
	\Omega_v(a,b) = \{ [x_0:x_1] \in \cP(w_0,w_1)(K_v) : [x_0:x_1]=[a:b]\}.
\end{equation}
Partition $\Omega_v^{\aff}(a,b)\cap \cO_{K,v}^2$ by the sets
\begin{equation}
\Omega_{v,t}^{\aff}(a,b)
 = \left\{(x_0,x_1)\in \cO_{K,v}^{2} : |x_0-q^{4t}a|_v\leq \frac{1}{q^{4t+1}},\ |x_1-q^{6t}b|_v\leq \frac{1}{q^{6t+1}}\right\}.
\end{equation}
Since $\Omega_{v,t}^{\aff}(a,b)$ have $v$-adic measure
\begin{equation}
	m_v\left(\Omega_{v,t}^{\aff}(a,b)\right)=\frac{1}{q_v^{10t+2}},
\end{equation}
the $v$-adic measure of $\Omega_{v}^{\aff}(a,b)\cap \cO_{K,v}^2$ is 
\begin{equation}
	m_v\left(\Omega_{v}^{\aff}(a,b)\cap \cO_{K,v}^2 \right) = \sum_{t\geq 0} \frac{1}{q^{10t+2}}=\frac{1}{q_v^2} \frac{q_v^{10}}{q_v^{10}-1}.
\end{equation} 

 Summing over the $3(q_v-1)/2$ pairs $(a,b)$ shows that the number of elliptic curves over $K$, with a point of order $4$, having $1$-multiplicative reduction with respect to $\phi$ at $v$, and of height bounded by $B$ is 
\begin{align}
\frac{q_v^{10}}{q_v^{10}-1} \kappa_{1}(G(1,4),v) \kappa_G B^{1/4}+O\left(B^{\frac{1}{4}-\frac{1}{12\deg_{\Q}(K)}}\right),
\end{align}
where $\kappa_{1}(G(1,4),v) = \frac{3q_v-3}{2q_v^2}$.
\end{example}


\section{Distribution of Selmer ratios}

In this section, we prove that the logarithms of ratios of sizes of certain Selmer groups are normally distributed (Theorem \ref{thm:selmer_erdos_kac_prime}), and use this to prove our main result (Theorem \ref{thm:main2}).

\subsection{Asymptotics for sums of Selmer ratios}\label{sec:Selmer_distribution}

Let $d=\deg(\phi)$, let $s(E)\defeq \log_d(\cS(E/E'))$, and let $s_v(E)\defeq \log_d(\cS_v(E/E'))$. 
Define
\begin{equation}
	s_G(B;K) \defeq \sum_{E\in \cE_G(B;K)} s(E) = \sum_{E\in \cE_G(B;K)} \sum_{v\in \Val(K)} s_v(E).
\end{equation}
 For 
\begin{equation}
	\star \in \{\textnormal{multiplicative (mult), additive (add), split multiplicative (split)}\}
\end{equation}
 let $s_{\star}$ denote the sum of $s_v(E)$ over $v$ for which $E$ has $\star$ reduction at $v$. Define
\begin{equation}
 s_{\star}(B;G,K) \defeq \sum_{E\in \cE_G(B;K)} s_{\star}(E).
\end{equation}


\begin{proposition}\label{prop:log_d_Tamagawa_estimate}
Suppose that the modular curve $\cX_G$ is isomorphic to a weighted projective stack, and let $\alpha(G)$ be as in Proposition \ref{prop:elliptic_count}. Then
\begin{equation}
s_G(B;K)=s_{mult}(B;G,K)+O_{K,G}\left(B^{\alpha(G)}\right).
\end{equation}
\end{proposition}

\begin{proof}
As $s_v(E)=0$ if $E$ has good reduction at $v$, we have 
\begin{equation}
s_G(B;K)=s_{mult}(B;G,K)+s_{add}(B;G,K).
\end{equation}
It therefore suffices to show $s_{add}(B;G,K)=O_{K,G}\left(B^{\alpha(T)}\right)$.

Let $S\subset \Val(K)$ be the set of all places above $2$, $3$, or primes dividing the level of $G$. For each elliptic curve $E$ define the set
\begin{equation}
S_{add}(E) \defeq \{v\in S : E \text{ has additive reduction at } v\}.
\end{equation}
By Proposition \ref{prop:local_selmer_ratio_bound} and Proposition \ref{prop:elliptic_count}, the contribution to $s_{add}(B;G,K)$ from places in $S$ is 
\begin{equation}\label{eq:bad_v_additive_estimate}
\sum_{E\in \cE_G(B;K)} \sum_{v\in S_{add}(E)} s_v(E)\leq \sum_{E\in\cE_G(B;K)} \sum_{v\in S_{add}(E)} 3 \ll_{K,G} \#\cE_G(B;K)\ll_{K,G}B^{\alpha(G)}.
\end{equation} 

 
By Proposition \ref{prop:elliptic_count}, the contribution to $s_{add}(B;G,K)$ from places not in $S$ is 
\begin{salign}\label{eq:good_v_additive_estimate}
s_{add}(B;G,K)-\sum_{E\in \cE_G(B;K)} \sum_{v\in S_{add}(E)} s_v(E)
& \ll_{K} \sum_{v\in \Val_0(K)} \kappa_{add}(G,\p) B^{\alpha(G)}\\
& \ll_{K,G} \sum_{v\in \Val_0(K)} \frac{B^{\alpha(G)}}{N_{K/\Q}(\p_v)^2}\\
& \ll_{K} B^{\alpha(G)}.
\end{salign}
The estimates (\ref{eq:bad_v_additive_estimate}) and (\ref{eq:good_v_additive_estimate}) together give the desired bound.
\end{proof}

\begin{corollary}
If $\deg(\phi)\neq 2$ then 
\begin{equation}
s_G(B;K)=s_{split}(B;G,K) + O_{K,G}(B^{\alpha(G)}).
\end{equation}
\end{corollary}

Suppose $\deg(\phi)=2$ and let $v\in \Val_0(K)$ be a place at which $E$ has multiplicative reduction. Define the \textbf{local pseudo Selmer ratio} to be
\begin{equation}\label{eq:local_pseudo_Selmer_ratio}
\tilde{\cS}_v(E/E')\defeq \frac{v(\Delta(E'))}{v(\Delta(E))}.
\end{equation}
Accordingly, define $\tilde{s}_v(E)\defeq\log_2(\tilde{\cS}_v(E/E')$, define $\tilde{s}_{\mult}(E)$ to be the sum of $\tilde{s}_v(E)$ over $v$ for which $E$ has multiplicative reduction, and define
\begin{equation}
\tilde{s}_{mult}(B;G,K) \defeq \sum_{E\in \cE_G(B;K)} \tilde{s}_{mult}(E). 
\end{equation}

\begin{corollary}
If $\deg(\phi)= 2$, then 
\begin{equation}
s_G(B;K)=\tilde{s}_{mult}(B;G,K) + O_{K,G}(B^{\alpha(G)}).
\end{equation}
\end{corollary}

\begin{proof}
By Proposition \ref{prop:tamagawa_ratio_values} the local pseudo Selmer ratio (\ref{eq:local_pseudo_Selmer_ratio}) equals the usual local Selmer ratio (\ref{eq:Selmer_ratio}) except when $v$ is above $2$ or when $2$ divides $\gcd(v(\Delta(E)),v(\Delta(E')))$; let $S(E)$ denote this set of places. These are the only places where the sums $s_{\mult}(E)$ and $\tilde{s}_{mult}(E)$ differ. Therefore, by Proposition \ref{prop:log_d_Tamagawa_estimate} it will suffice to show that
\begin{equation}\label{eq:pseudo_sum}
\sum_{E\in \cE_G(B;K)} \sum_{v\in S(E)} 1 \ll_G B^{\alpha(G)}.
\end{equation}

Let $S'(E)$ denote the subset of places of $S(E)$ above $2$, $3$, or primes dividing the level of $G$. The sum over $S'(E)$ can be dealt with as in the proof of Proposition \ref{prop:log_d_Tamagawa_estimate},
\begin{equation}\label{eq:psuedo_sum_estimate_2,3}
\sum_{E\in \cE_G(B;K)} \sum_{v\in S'(E)} s_v(E)\leq \sum_{E\in \cE_G(B;K)} \sum_{v\in S'(E)} 3 \ll_{K,G} \#\cE_G(B;K)\ll_{K,G} B^{\alpha(G)}.
\end{equation}

By \cite[Theorem 1.1.2]{Phi24+}, the proportion of elliptic curves $E\in \cE_G(K)$ with 
\begin{equation}
\p^2|\gcd(\Delta(E)),v(\Delta(E'))
\end{equation}
 for a fixed prime $\p\in \sP(K)$ not in $S'$ is 
\begin{salign}
&\frac{\#\{(A,B)\in \cO_K/\p^2\times \cO_K/\p^2 : f_{G,4}(A,B)\equiv f_{G,6}(A,B)\equiv 0\pmod{\p^2} \}}{\#\left(\cO_K/\p^2\times \cO_K/\p^2\right)}\\
&\qquad \ll \frac{N_{K/\Q}(\p)^2\max\{\deg(f_{G,4}),\deg(f_{G,6})\}}{N_{K/\Q}(\p)^4}\\
&\qquad \ll_G N_{K/\Q}(\p)^{-2}.
\end{salign}

This, together with Proposition \ref{prop:elliptic_count}, shows that the contribution to the sum (\ref{eq:pseudo_sum}) from places not in $S'$ is
\begin{salign}\label{eq:psuedo_sum_estimate}
\sum_{E\in \cE_G(B;K)}\sum_{v\in S(E)} s_v(E)-\sum_{E\in \cE_G(B;K)}\sum_{v\in S_{2,3}(E)} s_v(E)
& \ll_{K,G} \sum_{v\in \Val_0(K)} \frac{B^{\alpha(G)}}{N_{K/\Q}(\p_v)^2}\\
&\ll_{K} B^{\alpha(T)}.
\end{salign}
The estimates (\ref{eq:psuedo_sum_estimate_2,3}) and (\ref{eq:psuedo_sum_estimate}) together give the desired bound (\ref{eq:pseudo_sum}).
\end{proof}


\subsection{Isogenies of prime degree}

In this subsection we focus on the case that the isogeny $\phi$ is of prime degree. For some examples of isogenies of composite degree see Section \ref{sec:Composite_case}.

\begin{theorem}\label{thm:selmer_erdos_kac_prime}
	Consider the family of elliptic curves $E$ over a number field $K$ with level structure $G=G(M,MN)$. Let $\Phi$ be a subgroup of $\Z/M\Z$ of prime order, and let $m$ be as in Table \ref{tab:mean_variance_II}. Let $\sG = \sG_m(K)$ denote the subgroup 
\begin{equation}
	\{a\in (\Z/m\Z)^\times : \exists \text{ infinitely many degree one primes $\p$ with } N_{K/\Q}(\p)\equiv a \pmod{m}\}.
\end{equation}
	Let $c_{\E}(G,\Phi,\sG)$ and $\c_{\V}(G,\Phi,sG)$ be the constants in Table \ref{tab:mean_variance_II}. 
	Then the logarithmic Selmer ratios $s(E)$ are normally distributed with mean
	\begin{equation}
		c_{\E}(G,\Phi,\sG)\log\log(B)+O(1)
	\end{equation}
	and variance
	\begin{equation}
		c_{\V}(G,\Phi,\sG)\log\log(B)+O(1).
	\end{equation}
\end{theorem}

\begin{proof}
To study the distribution of the logarithmic Selmer ratios $s(E)$ we follow the ideas of Billingsley's proof of the Erd\H{o}s--Kac Theorem \cite{Bil69}. The idea is that we would like to apply the central limit theorem to the logarithmic Selmer ratios $s(E)$, but the hypotheses are not met. So we instead model the situation with independent random variables for which the central limit theorem can be applied directly. Finally, we use the method of moments to show that the distribution of the $s(E)$ agrees with the distribution coming from the random variables.

Fix a level structure $G=G(M,MN)\in \bT$ and fix a subgroup $\Phi$ of $\Z/M\Z$ of prime order. Let $K$ be a number field for which the modular curve $\cX_G$ has a $K$-rational point.
We model the distribution of
\begin{equation}
	s_+(B) \defeq \sum_{N_{K/\Q}(\p_v)\leq B} \#\{E\in \cE(K; B) : s_v(E)=1\}
\end{equation}
and
\begin{equation}
	s_-(B) \defeq \sum_{N_{K/\Q}(\p_v)\leq B} \#\{E\in \cE(K; B) : s_v(E)=-1\}
\end{equation}
using random variables.

  Let $v\in \Val_0(K)$ be a finite place of $K$. If $v$ is above $2$, $3$, $\deg(\phi)$, or any prime dividing $MN$, then define $\kappa_{-}(G,\Phi,v)=0$ and $\kappa_{+}(G,\Phi,v)=0$. Otherwise, let $\kappa_{-}(G,\Phi,v)$ and $\kappa_{+}(G,\Phi,v)$ be as in Table \ref{tab:isogeny_local_densities}. Define independent Bernoulli random variables
\begin{equation}
X_{+,v}\defeq \begin{cases}
1 & \text{ with probability } \kappa_+(G,\Phi,v),\\
0 & \text{ with probability } 1-\kappa_+(G,\Phi,v),
\end{cases}
\end{equation}
and
\begin{equation}
X_{-,v}\defeq \begin{cases}
1 & \text{ with probability } \kappa_-(G,\Phi,v),\\
0 & \text{ with probability } 1-\kappa_-(G,\Phi,v).
\end{cases}
\end{equation}
 Let $\varepsilon(E,a_v,m)$ is as in Proposition \ref{prop:K-Dirichlet}, and define
\begin{equation}
X_+(B)\defeq \sum_{N_{K/\Q}(\p_v)\leq B} \varepsilon(K,a_v,m) X_{+,v}
\end{equation}
and 
\begin{equation}
X_-(B)\defeq \sum_{N_{K/\Q}(\p_v)\leq B} \varepsilon(K,a_v,m) X_{-,v}.
\end{equation}
Proposition \ref{prop:log_d_Tamagawa_estimate} gives justification for why these random variables should model $s_+(B)$ and $s_-(B)$ despite not taking primes of additive reduction into account. 

Applying the multidimensional central limit theorem to the random vectors $\begin{bmatrix} X_{+,v} & X_{-,v} \end{bmatrix}
$ shows that $X_+(B)$ and $X_-(B)$ converge to independent normal distributions.
Let $c_+(G,\Phi,\sG)$ and $c_-(G,\Phi,\sG)$ be as in Table \ref{tab:mean_variance_I}. Then the distribution of $X_+(B)$ has mean and variance 
\begin{equation}
    c_+(G,\Phi,\sG)\log\log(B),
\end{equation}
and the distribution for $X_-(B)$ has mean and variance 
\begin{equation}
    c_-(G,\Phi,\sG)\log\log(B).
\end{equation}  

We now show that the distributions described by $X_+(B)$ and $X_-(B)$ coincide with the distributions described by $s_+(B)$ and $s_-(B)$. For this we use the method of moments, i.e., we show that all of the mixed moments agree.

Let $\Delta_+(E)$ denote the largest factor of $\Delta(E)$ for which $v(\Delta(\varphi(E)))=\deg(\phi) v(\Delta(E))$, and let $\Delta_-(E)$ denote the largest factor of $\Delta(E)$ for which $\deg(\phi) v(\Delta(\varphi(E)))=v(\Delta(E))$.\footnote{See Table \ref{tab:Discriminants} for explicit expressions for $\Delta(E)$ in the cases of interest.} Then
\begin{equation}
	s_\pm(B) = \sum_{N_{K/\Q}(\p_v)\leq B} \#\{E\in \cE(K; B) : \p_i|\Delta_\pm(E)\}.
\end{equation}

By the definition of $X_+(B)$ and $X_-(B)$, for any $m,n\in \Z_{\geq 0}$ we have
\begin{align*}
\E\left[s_+(B)^{m} s_-(B)^{n}\right]
&= \lim_{B\to\infty} \sum_{\substack{\p_1,\dots,\p_{m}\in \sP(B;K)\\ \q_1,\dots,\q_{n}\in\sP(B;K)}} \frac{\#\{E\in \cE_G(B;K) : \p_i|\Delta_+(E) \text{ and } \q_j|\Delta_-(E)\ \forall i,j\}}{\# \cE_G(B;K)}\\
&= \E\left[X_+(B)^{m} X_-(B)^{n}\right]+O(1).
\end{align*}
Similarly,
\begin{salign}
&\E\left[ (s_+(M)-\E [s_+(B)])^{m}(s_-(B)-\E(s_-(B)))^{n}\right]\\
&\qquad =\E\left[(X_+(B)-\E[X_+(B)])^{m}(X_-(B)-\E[X(B)])^{n}\right]+O(1).
\end{salign}

As the difference of normal distributions is still a normal distribution, we conclude that $s(B)\sim s_+(B)-s_-(B)$ tends to a normal distribution. 
Let $c_\E(G,\Phi,\sG)$ and $c_\V(G,\Phi,\sG)$  be as in Table \ref{tab:mean_variance_II}. 
 Then, as $B\to\infty$, the $s(B)$ are normally distributed with mean $c_\E(G,\Phi,K)\log\log(B)+O(1)$ and variance $c_\V(G,\Phi,K)\log\log(B)+O(1)$.
\end{proof}

We now prove our main result, Theorem \ref{thm:main}, which we restate for the reader's convenience. 

\begin{theorem}\label{thm:main2}
Suppose that $X_1(M,MN)$ has genus zero, and suppose that $p$ is a prime divisor of $MN$. Let $\theta(G, \Phi, \sG_m(K))$ denote the constant in Table \ref{tab:mean_variance_II}. Then
\begin{equation}\label{eq:main2}
	\limsup_{B\to\infty} \frac{1}{\# \cE_{G}(B; K)} \sum_{E\in \cE_{G}(B;K)} \# \Sel_p(E) \gg \log(B)^{\theta(G, \Phi, \sG_m(K))}.
\end{equation}
\end{theorem}

\begin{proof}
	By Theorem \ref{thm:selmer_erdos_kac_prime} we know that as $B\to\infty$ the logarithmic Selmer ratios are normally distributed with mean 
    \begin{equation}
        c_{\E}(G,\Phi, \sG_m(K))\log\log(B)+O(1)
    \end{equation} 
    and variance 
    \begin{equation}
        c_{\V}(G,\Phi, \sG_m(K))\log\log(B)+O(1).
    \end{equation} 
    It follows that as $B\to\infty$ the Selmer ratios $\#\Sel_{\phi}(E)/\#(\Sel_{\widehat{\phi}}(E')$ become lognormally distributed with mean
	\begin{salign}
		&\exp\left(c_{\E}(G,\Phi, \sG_m(K))\log\log(B)+c_{\V}(G,\Phi, \sG_m(K))\log\log(B)/2 +O(1)\right) \\
		&= \exp\left(\left(c_{\E}(G,\Phi,\sG_m(K)) + \frac{c_{\V}(G,\Phi, \sG_m(K))}{2}\right)\log\log(B) + O(1)\right)\\
 &\asymp  \log(B)^{c_{\E}(G,\Phi,\sG_m(K)) + c_{\V}(G,\Phi, \sG_m(K))/2}. 
	\end{salign}
	By Equation (\ref{eq:N-Selmer_inequality}) we obtain the desired result upon 
	setting
	\begin{equation}	 
	\theta(G,\Phi,\sG_m(K)) = c_{\E}(G,\Phi,\sG_m(K)) + c_{\V}(G,\Phi, \sG_m(K))/2.
	\end{equation}
\end{proof}

\section{Isogenies of composite degree}\label{sec:Composite_case}

The methods used to prove Theorem \ref{thm:selmer_erdos_kac_prime} also work in the case of isogenies of composite degree. In the case of prime isogenies, the local, logarithmic, Selmer ratios are all $-1$, $0$, or $1$. However, for isogenies of composite degree there are more possibilities. As an example, for isogenies of degree $9$ the possible local, logarithmic, Selmer ratios are 
\begin{center}
$-1$, $-1/2$, $0$, $1/2$, and $1$.
\end{center}
 This becomes even more cumbersome for isogenies that are not of prime power degree. As an example, for isogenies of degree $6$ the possible local, logarithmic, Selmer ratios are 
 \begin{center}
 $-1$, $-\log_6(3)$, $-\log_6(2)$, $-\log_6(3/2)$, $0$, $\log_6(3/2)$, $\log_6(2)$,  $\log_6(3)$, and $1$.
 \end{center}
  Because of this, the results become more cluttered. In this section we give only a couple examples in this direction as an illustration (although one could easily obtain results in more cases if they so desired). 

\subsection{9-Selmer groups of elliptic curves with 9-torsion}\label{subsec:9-selmer}

Let $E$ be an elliptic curve with a $K$ rational $9$-torsion point $P\in E[9](K)$, and let $\phi: E\to E'$ be the isogeny of degree $9$ with kernel $\langle P\rangle$. 

In Table \ref{tab:9-isogeny_local_densities} we give densities for elliptic curves over $K$ with $G(1,9)$-level structure and with $n$-reduction at $\phi$, for each $n\in \{-1, -1/2, 0, 1/2, 1\}$. These densities are computed as in Section \ref{sec:local_conditions}. 

\begin{longtblr}[caption = {Local densities of elliptic curves with $G(1,9)$-level structure over $K$ with degree $9$ isogeny $\phi$ at places not above $2$ or $3$.},
label = {tab:9-isogeny_local_densities}]{
  width = \linewidth, rowhead = 2,
  row{1} = {white}, 
colspec={cc},}
\toprule[1.5pt]
Conditions on & $n$-reduction at $\phi$  \\
 $N_{K/\Q}(\p_v)=q$ & $\left(\kappa_{-1}(G,\phi_\Phi, v), \kappa_{-\frac{1}{2}}(G,\phi_\Phi, v), \kappa_{0}(G,\phi_\Phi, v), \kappa_{\frac{1}{2}}(G,\phi_\Phi, v), \kappa_{1}(G,\phi_\Phi, v)\right)$ \\
\midrule
%
%
%
%
 $q\equiv 1 \pmod{9}$  & $\left(\frac{3q-3}{q^2},0,\frac{2q-2}{q^2},0,\frac{3q-3}{q^2}\right)$ \\
 $q\equiv 2, 5 \pmod{9}$  & $\left(\frac{3q-3}{q^2},0,0,0,0\right)$ \\
 $q\equiv  4, 7 \pmod{9}$  &  $\left(\frac{3q-3}{q^2},0,\frac{2q-2}{q^2},0,0\right)$ \\
 $q\equiv 8 \pmod{9}$  &  $\left(\frac{3q-3}{q^2},0,\frac{3q-3}{q^2},0,0\right)$  \\
\bottomrule[1.5pt]
\end{longtblr}

The method used to prove Theorem \ref{thm:selmer_erdos_kac_prime} works in this case as well. 
 For each $n\in \{-1,-1/2,0,1/2,1\}$ define a random variable
 \begin{equation}
X_{n,v}\defeq \begin{cases}
1 & \text{ with probability } \kappa_n(G,\Phi,v),\\
0 & \text{ with probability } 1-\kappa_n(G,\Phi,v).
\end{cases}
\end{equation}
Then we model the distribution of 
\begin{equation}
	s_n(B) \defeq \sum_{N_{K/\Q}(\p_v)\leq B} \#\{E\in \cE(K; B) : s_v(E)=n\}
\end{equation}
by the sum of random variables
\begin{equation}
X_n(B)\defeq \sum_{N_{K/\Q}(\p_v)\leq B} \varepsilon(K,a_v,m) X_{n,v}.
\end{equation}
Then we apply the Central Limit Theorem to obtain the distribution of the $X_n(B)$. We note the distribution will depend on the subgroup of $9$-th roots of unity contained in $K$, which we denote by $\sG=\sG_9(K)$. Table \ref{tab:9-isogeny_distributions} gives constants $c_n(G(1,9),C_9,\sG)$, such that $X_n(B)$ has mean and variance
\begin{equation}
	c_n\left(G(1,9),C_9,\sG\right) \log\log(B) + O(1).
\end{equation}

  Using the method of moments one can show that the distributions of $s_n$ and $X_n$ are the same.

\begin{longtblr}[caption = {Mean and variance constants for $X_n(B)$.},
label = {tab:9-isogeny_distributions}]{
  width = \linewidth, rowhead = 2,
  row{1} = {white}, 
colspec={cc},}
\toprule[1.5pt]
Subgroup & Mean and variance constants \\
 $\sG_9(K)$ & $\left(c_{-1}(G(1,9),C_9,\sG_9(K)), \dots, c_1(G(1,9),C_9,\sG_9(K))\right)$ \\
\midrule
%
%
%
%
 $\{1\}$  & $\left(3,0,2,0,3\right)$ \\
 $\{1,8\}$  & $\left(3,0,\frac{5}{2}, 0, \frac{3}{2}\right)$ \\
 $\{1,4,7\}$  &  $\left(3,0,3,0,1\right)$ \\
 $(\Z/9\Z)^\times$  &  $\left(3,0,\frac{3}{2},0,\frac{1}{2}\right)$  \\
\bottomrule[1.5pt]
\end{longtblr}

These can be combined to give the distribution of the logarithmic Selmer ratios. Let $c_{\E}(G(1,9),C_9,\sG_9(K))$ and $c_{\V}(G(1,9),C_9,\sG_9(K))$ be as in Table \ref{tab:9-isogeny_selmer_distributions}. Then, as $B\to\infty$, the logarithmic Selmer ratios 
\begin{equation}
	s(B) = \sum s_n(B)
\end{equation}
 are normally distributed with mean 
\begin{equation}
	c_{\E}(G(1,9),C_9,\sG_9(K))\log\log(B) + O(1)
\end{equation}
and variance
\begin{equation}
	c_{\V}(G(1,9),C_9,\sG_9(K))\log\log(B) + O(1).
\end{equation}

\begin{longtblr}[ caption = {Mean and variance constants for $s(B)$.},
label = {tab:9-isogeny_selmer_distributions}]{
  width = \linewidth, rowhead = 2,
  row{1} = {white}, 
colspec={cccc},}
\toprule[1.5pt]
Subgroup & Mean constant & Variance constant &  \\
 $\sG_9(K)$ & $c_{\E}(G(1,9),C_9,\sG_9(K))$ & $c_{\V}(G(1,9),C_9,\sG_9(K))$ &  $\theta(G(1,9),C_9,\sG_9(K))$ \\
\midrule
%
%
%
%
 $\{1\}$  & $0$ & $6$ & $3$ \\
 $\{1,8\}$  & $-\frac{3}{2}$ & $\frac{9}{2}$ & $\frac{3}{4}$ \\
 $\{1,4,7\}$  & $-2$ & $4$ & $0$ \\
 $(\Z/9\Z)^\times$  & $-\frac{5}{2}$ & $\frac{7}{2}$ & $- \frac{3}{4}$ \\
\bottomrule[1.5pt]
\end{longtblr}

Let $\theta(G(1,9),C_9,\sG_9(K))$ be as in Table \ref{tab:9-isogeny_selmer_distributions}. 
By the same argument as the proof of Theorem \ref{thm:main2}, we have
\begin{equation}
	\limsup_{B\to\infty} \frac{1}{\# \cE_{G(1,9)}(B; K)} \sum_{E\in \cE_{G(1,9)}(B;K)} \# \Sel_9(E) \gg \log(B)^{\theta(G(1,9), C_9, \sG_9(K))}.
\end{equation}

\subsection{6-Selmer groups of elliptic curves with 6-torsion}\label{subsec:6-selmer}

Now let $E$ be an elliptic curve with a $K$ rational $6$-torsion point $P\in E[6](K)$, and let $\phi: E\to E'$ be the degree $6$ isogeny with kernel $\langle P\rangle$. 

In Table \ref{tab:6-isogeny_local_densities} we give densities for elliptic curves over $K$ with $G(1,9)$-level structure and with $n$-reduction at $\phi$, for each 
\begin{equation}\label{eq:n_6-selmer}
n\in \{-1, -\log_6(3), -\log_6(2),  -\log_6(3/2), 0, \log_6(3/2), \log_6(2), \log_6(3), 1\}.
\end{equation}
 These densities are computed as in Section \ref{sec:local_conditions}. 

\begin{longtblr}[ caption = {Local densities of elliptic curves with $G(1,6)$-level structure over $K$ with degree $6$ isogeny $\phi$ at places not above $2$ or $3$.},
label = {tab:6-isogeny_local_densities}]{
  width = \linewidth, rowhead = 2,
  row{1} = {white}, 
colspec={cc},}
\toprule[1.5pt]
Conditions on & $n$-reduction at $\phi$  \\
 $N_{K/\Q}(\p_v)=q$ & $\left(\kappa_{-1}(G,\phi_\Phi, v), \dots , \kappa_{1}(G,\phi_\Phi, v)\right)$ \\
\midrule
%
%
%
%
 $q\equiv 1 \pmod{6}$  & $\left(\frac{q-1}{q^2},0,0,\frac{q-1}{q^2},0,\frac{q-1}{q^2},0,0,\frac{q-1}{q^2}\right)$ \\
 $q\equiv 5 \pmod{6}$  & $\left(\frac{q-1}{q^2},0,\frac{q-1}{q^2},\frac{q-1}{q^2},0,0,\frac{q-1}{q^2},0,0\right)$ \\
\bottomrule[1.5pt]
\end{longtblr}

The method used to prove Theorem \ref{thm:selmer_erdos_kac_prime} works in this case as well. For $n$ as in (\ref{eq:n_6-selmer}), define $X_{n}(B)$ and $s_n(B)$ as in Section \ref{subsec:9-selmer}. Let $c_n(G(1,6),C_6,\sG_6(K))$ be the as in Table \ref{tab:6-isogeny_distributions}. An argument using the central limit theorem shows that $X_n(B)$ becomes normally distributed with mean and variance
\begin{equation}
	c_n(G(1,6),C_6,\sG_6(K)) \log\log(B) + O(1).
\end{equation}

\begin{longtblr}[ caption = {Mean and variance constants for $X_n(B)$.},
label = {tab:6-isogeny_distributions}]{
  width = \linewidth, rowhead = 2,
  row{1} = {white}, 
colspec={cc},}
\toprule[1.5pt]
Subgroup & Mean and variance constants \\
 $\sG(K,6)$ & $\left(c_{-1}(G(1,6),C_6,\sG_6(K)), \dots, c_1(G(1,6),C_6,\sG_6(K))\right)$ \\
\midrule
%
%
%
%
 $\{1\}$  & $\left(1,0,0,1,0,1,0,0,1\right)$ \\
 $(\Z/6\Z)^\times$  &  $\left(1,0,\frac{1}{2},1,0,\frac{1}{2},\frac{1}{2},0,\frac{1}{2}\right)$  \\
\bottomrule[1.5pt]
\end{longtblr}

Using the method of moments one can show that the distributions of $s_n$ and $X_n$ are the same.

Let $c_{\E}(G(1,6),C_6,\sG_6(K))$ and $c_{\V}(G(1,6),C_6,\sG_6(K))$ be as in Table \ref{tab:6-isogeny_selmer_distributions}. Then, as $B\to\infty$, the logarithmic Selmer ratios 
\begin{equation}
	s(B) = \sum s_n(B)
\end{equation}
 are normally distributed with mean 
\begin{equation}
	c_{\E}(G(1,6),C_6,\sG_6(K))\log\log(B) + O(1)
\end{equation}
and variance
\begin{equation}
	c_{\V}(G(1,6),C_6,\sG_6(K))\log\log(B) + O(1).
\end{equation}

\begin{longtblr}[ caption = {Mean and variance constants for $s(B)$.},
label = {tab:6-isogeny_selmer_distributions}]{
  width = \linewidth, rowhead = 2,
  row{1} = {white}, 
colspec={ccc},}
\toprule[1.5pt]
Subgroup & Mean constant & Variance constant \\
 $\sG(K,6)$ & $c_{\E}(G(1,6),C_6,\sG_6(K))$ & $c_{\V}(G(1,6),C_6,\sG_6(K))$ \\
\midrule
%
%
%
%
 $\{1\}$  & $0$ & $2\left(1+\log_6\left(\frac{3}{2}\right)^2\right)$ \\
 $(\Z/3\Z)^\times$  & $-\frac{1}{2} \left(1 + \log_6\left(\frac{3}{2}\right) \right)$ & $\frac{1}{2}\left(3+2\log_6(2)^2+3\log_6\left(\frac{3}{2}\right)^2\right)$  \\
\bottomrule[1.5pt]
\end{longtblr}

Define
\begin{align}
	\theta(G(1,6),C_6,\{1\})&\defeq 1+\log_6\left(\frac{3}{2}\right)^2\\
	\theta(G(1,6),C_6,\{1\})&\defeq \frac{1}{4} - \frac{1}{2}\log_6\left(\frac{3}{2}\right) + \frac{1}{2} \log_6(2)^2 + \frac{3}{4} \log_6\left(\frac{3}{2}\right)^2.
\end{align}
By the same argument as the proof of Theorem \ref{thm:main2}, we have
\begin{equation}
	\limsup_{B\to\infty} \frac{1}{\# \cE_{G(1,6)}(B; K)} \sum_{E\in \cE_{G(1,6)}(B;K)} \# \Sel_6(E) \gg \log(B)^{\theta(G(1,6), C_6, \sG_6(K))}.
\end{equation}

\section{Tables}\label{sec:tables}

\begin{longtblr}[ caption = {Defining equations for torsion families in short Weierstrass form.},label = {tab:Weierstrass_families},]{
  width = \linewidth, stretch=1.5, colspec = {ccX[c]}, rowhead = 1, row{2-3,6-7, 10-11, 14-15, 18-19, 22-23, 26-27, 30-31, 34-35} = {gray!20}, row{1} = {white},}
\toprule[1.5pt] 
$G$ &  $i$ & $f_{G,i}(A,B)$  \\
\midrule
$G(1,2)$ &   $4$ &  $-432(A^2 - 3B)$\\
			  &      $6$ & $1728A(2A^2 - 9B)$\\ 
$G(1,3)$ &   $4$ &  $-27A(A^3 - 24B)$\\
 			  &      $6$ & $54(A^6 - 36A^3B + 216B^2)$\\ 
$G(1,4)$ &   $4$ &  $-27(B^{4} + 16 A B^{2} + 16 A^{2})
$\\
  &      $6$ & $54(B^{2} + 8 A)   (B^{4} + 16 A B^{2} - 8 A^{2})
$\\ 
$G(1,5)$ &   $4$ &  $-27(A^4 - 12A^3B + 14A^2B^2 + 12AB^3 + B^4)$\\
 &      $6$ & $54(A^2+B^2)(A^4 - 18A^3B + 74A^2B^2 + 18AB^3 + B^4)$\\ 
$G(1,6)$ &   $4$ &  $-27(3A + B)  (3A^3 + 3A^2B + 9AB^2 + B^3)
$\\
  &      $6$ & $-54(3A^2 - 6AB - B^2)  (9A^4 + 36A^3B + 30A^2B^2 + 12AB^3 + B^4)$\\ 
%
%
%
%
$G(1,7)$ &   $4$ &  $-27(A^2 - AB + B^2)   (A^6 - 11A^5B + 30A^4B^2 - 15A^3B^3 - 10A^2B^4 + 5AB^5 + B^6)$\\
%
  &      $6$  & $54(A^{12} - 18A^{11}B + 117A^{10}B^2 - 354A^9B^3 + 570A^8B^4 - 486A^7B^5 + 273A^6B^6 - 222A^5B^7 + 174A^4B^8 - 46A^3B^9 - 15A^2B^{10} + 6AB^{11} + B^{12})$\\ 
%
%
%
$G(1,8)$ &   $4$ &  $-27(16 A^{8} - 64 A^{7} B + 224 A^{6} B^{2} - 448 A^{5} B^{3} + 480 A^{4} B^{4} - 288 A^{3} B^{5} + 96 A^{2} B^{6} - 16 A B^{7} + B^{8})
$\\
 &      $6$ & $-54(8 A^{4} - 16 A^{3} B + 16 A^{2} B^{2} - 8 A B^{3} + B^{4})  (8 A^{8} - 32 A^{7} B - 80 A^{6} B^{2} + 352 A^{5} B^{3} - 456 A^{4} B^{4} + 288 A^{3} B^{5} - 96 A^{2} B^{6} + 16 A B^{7} - B^{8})
$\\ 
%
%
%
$G(1,9)$ &   $4$ &  $-27(A^{3} - 3 A^{2} B + B^{3})   (A^{9} - 9 A^{8} B + 27 A^{7} B^{2} - 48 A^{6} B^{3} + 54 A^{5} B^{4} - 45 A^{4} B^{5} + 27 A^{3} B^{6} - 9 A^{2} B^{7} + B^{9})
$\\
 &      $6$ & $54(A^{18} - 18 A^{17} B + 135 A^{16} B^{2} - 570 A^{15} B^{3} + 1557 A^{14} B^{4} - 2970 A^{13} B^{5} + 4128 A^{12} B^{6} - 4230 A^{11} B^{7} + 3240 A^{10} B^{8} - 2032 A^{9} B^{9} + 1359 A^{8} B^{10} - 1080 A^{7} B^{11} + 735 A^{6} B^{12} - 306 A^{5} B^{13} + 27 A^{4} B^{14} + 42 A^{3} B^{15} - 18 A^{2} B^{16} + B^{18})
$\\ 
%
%
%
%
$G(1,10)$ &   $4$ &  $-27(16 A^{12} - 128 A^{11} B + 416 A^{10} B^{2} - 720 A^{9} B^{3} + 720 A^{8} B^{4} - 288 A^{7} B^{5} - 256 A^{6} B^{6} + 432 A^{5} B^{7} - 240 A^{4} B^{8} + 40 A^{3} B^{9} + 16 A^{2} B^{10} - 8 A B^{11} + B^{12})
$\\
 &      $6$ & $54(2 A^{2} - 2 A B + B^{2})   (2 A^{4} - 2 A B^{3} + B^{4})   (4 A^{4} - 12 A^{3} B + 6 A^{2} B^{2} + 2 A B^{3} - B^{4})   (4 A^{8} - 32 A^{7} B + 104 A^{6} B^{2} - 176 A^{5} B^{3} + 146 A^{4} B^{4} - 48 A^{3} B^{5} - 4 A^{2} B^{6} + 6 A B^{7} - B^{8})$\\ 
%
%
%
%
%
$G(1,12)$ &   $4$ &  $-27(6 A^{4} - 12 A^{3} B + 12 A^{2} B^{2} - 6 A B^{3} + B^{4})   (24 A^{12} - 144 A^{11} B + 864 A^{10} B^{2} - 3000 A^{9} B^{3} + 6132 A^{8} B^{4} - 8112 A^{7} B^{5} + 7368 A^{6} B^{6} - 4728 A^{5} B^{7} + 2154 A^{4} B^{8} - 684 A^{3} B^{9} + 144 A^{2} B^{10} - 18 A B^{11} + B^{12})$\\
 &      $6$ & $-54(24 A^{8} - 96 A^{7} B + 216 A^{6} B^{2} - 312 A^{5} B^{3} + 288 A^{4} B^{4} - 168 A^{3} B^{5} + 60 A^{2} B^{6} - 12 A B^{7} + B^{8})   (72 A^{16} - 576 A^{15} B - 1008 A^{14} B^{2} + 17136 A^{13} B^{3} - 65880 A^{12} B^{4} + 146304 A^{11} B^{5} - 222552 A^{10} B^{6} + 248688 A^{9} B^{7} - 211296 A^{8} B^{8} + 138720 A^{7} B^{9} - 70632 A^{6} B^{10} + 27696 A^{5} B^{11} - 8208 A^{4} B^{12} + 1776 A^{3} B^{13} - 264 A^{2} B^{14} + 24 A B^{15} - B^{16})$\\ 
%
%
%
%
$G(2,2)$ &   $4$ &  $-432(A^{2} - A B + B^{2})
$\\
 &      $6$ & $-1728(A - 2 B)   (A + B)   (2 A - B)
$\\ 
%
%
%
%
$G(2,4)$ &   $4$ &  $-432(256 A^{4} + 224 A^{2} B^{2} + B^{4})
$\\
 &      $6$ & $-3456(16 A^{2} + B^{2})   (16 A^{2} - 24 A B + B^{2})   (16 A^{2} + 24 A B + B^{2})
$\\ 
%
%
%
%
$G(2,6)$ &   $4$ &  $-432(3 A^{2} + B^{2})   (3 A^{6} + 75 A^{4} B^{2} - 15 A^{2} B^{4} + B^{6})
$\\
 &  $6$ & $-3456(3 A^{4} + 6 A^{2} B^{2} - B^{4})   (3 A^{4} - 24 A^{3} B + 6 A^{2} B^{2} - B^{4})   (3 A^{4} + 24 A^{3} B + 6 A^{2} B^{2} - B^{4})
$\\ 
%
%
%
%
$G(2,8)$ &   $4$ &  $-432(A^{16} - 8 A^{14} B^{2} + 12 A^{12} B^{4} + 8 A^{10} B^{6} + 230 A^{8} B^{8} + 8 A^{6} B^{10} + 12 A^{4} B^{12} - 8 A^{2} B^{14} + B^{16})
$\\
 &      $6$ & $3456(A^{8} - 4 A^{6} B^{2} - 26 A^{4} B^{4} - 4 A^{2} B^{6} + B^{8})   (A^{8} - 4 A^{6} B^{2} - 2 A^{4} B^{4} - 4 A^{2} B^{6} + B^{8})   (A^{8} - 4 A^{6} B^{2} + 22 A^{4} B^{4} - 4 A^{2} B^{6} + B^{8})
$\\ 
%
%
%
%
$G(3,3)$ &   $4$ &  $-27A   (A + 6 B)   (A^{2} - 6 A B + 36 B^{2})
$\\
 &      $6$ & $54(A^{2} - 6 A B - 18 B^{2})   (A^{4} + 6 A^{3} B + 54 A^{2} B^{2} - 108 A B^{3} + 324 B^{4})
$\\ 
%
%
%
%
$G(3,6)$ &   $4$ &  $-27(81 A^{3} + 27 A^{2} B + 9 A B^{2} + B^{3})   (243 A^{3} + 81 A^{2} B + 9 A B^{2} + B^{3})   (19683 A^{6} + 13122 A^{5} B + 2187 A^{4} B^{2} + 81 A^{2} B^{4} + 18 A B^{5} + B^{6})
$\\
 &      $6$ & $-54(27 A^{2} - B^{2})   (729 A^{4} + 486 A^{3} B + 108 A^{2} B^{2} + 18 A B^{3} + B^{4})   (729 A^{4} + 486 A^{3} B + 162 A^{2} B^{2} + 18 A B^{3} + B^{4})   (531441 A^{8} + 354294 A^{7} B + 118098 A^{6} B^{2} + 39366 A^{5} B^{3} + 10206 A^{4} B^{4} + 1458 A^{3} B^{5} + 162 A^{2} B^{6} + 18 A B^{7} + B^{8})
$\\ 
 %
%
%
%
$G(4,4)$ &   $4$ &  $-432(4 A^{4} - 12 A^{3} B + 14 A^{2} B^{2} - 6 A B^{3} + B^{4})   (4 A^{4} - 4 A^{3} B + 2 A^{2} B^{2} - 2 A B^{3} + B^{4})
$\\
 &      $6$ & $-3456\left(64\right)   (2 A^{2} - B^{2})   (2 A^{2} - 4 A B + B^{2})   (2 A^{4} - 4 A^{3} B + 6 A^{2} B^{2} - 4 A B^{3} + B^{4})   (8 A^{4} - 16 A^{3} B + 12 A^{2} B^{2} - 4 A B^{3} + B^{4})
$\\ 
%
%
%
%
$G(5,5)$ &   $4$ &  $-27(A^{4} + 3 A^{3} B - A^{2} B^{2} - 3 A B^{3} + B^{4})   (A^{8} - 4 A^{7} B + 7 A^{6} B^{2} - 2 A^{5} B^{3} + 15 A^{4} B^{4} + 2 A^{3} B^{5} + 7 A^{2} B^{6} + 4 A B^{7} + B^{8})   (A^{8} + A^{7} B + 7 A^{6} B^{2} - 7 A^{5} B^{3} + 7 A^{3} B^{5} + 7 A^{2} B^{6} - A B^{7} + B^{8})
$\\
 &      $6$ & $54(A^{2} + B^{2})   (A^{4} - 2 A^{3} B - 6 A^{2} B^{2} + 2 A B^{3} + B^{4})   (A^{8} - A^{6} B^{2} + A^{4} B^{4} - A^{2} B^{6} + B^{8})   (A^{8} - 4 A^{7} B + 17 A^{6} B^{2} - 22 A^{5} B^{3} + 5 A^{4} B^{4} + 22 A^{3} B^{5} + 17 A^{2} B^{6} + 4 A B^{7} + B^{8})   (A^{8} + 6 A^{7} B + 17 A^{6} B^{2} + 18 A^{5} B^{3} + 25 A^{4} B^{4} - 18 A^{3} B^{5} + 17 A^{2} B^{6} - 6 A B^{7} + B^{8})
$\\ 
\bottomrule[1.5pt]
\end{longtblr}


\begin{longtblr}[
   caption = {Discriminants for torsion families.}, 
   label = {tab:Discriminants}
]{
   width = \textwidth,  
   rowhead = 1,
   cells   = {font = \fontsize{10pt}{11pt}\selectfont},
   row{2-3,6-7, 10-12, 15-16, 19-21, 25-26, 29-31, 34-35, 39-40} = {gray!20},
   row{1} = {white},
   colspec = {c c X[c]}
}
\toprule[1.5pt]
$G$ & $d$ & Discriminant  \\
\midrule
$G(1,2)$ & $1$ & $2^4 B^2 (A^2-4B)$\\
   & $2$ & $2^8 B (A^2-4B)^2$\\ 
$G(1,3)$ & $1$ & $B^3(A^3-27B)$\\
   & $3$ & $B(A^3-27B)^3$\\ 
$G(1,4)$ & $1$ & $A^4B^2(16A+B^2)$\\
   & $2$ & $A^2B^4(16A+B^2)^2$\\
$G(1,5)$ & $1$ & $A^5B^5(A^2-11AB-B^2)$\\
   & $5$ & $AB(A^2-11AB-B^2)^5$\\ 
$G(1,6)$ & $1$ & $A^6B^2(A+B)^3(9A+B)$\\
			  & $2$ & $-A^3B(A+B)^6(9A+B)^2$\\
			  & $3$ & $A^2B^6(A+B)(9A+B)^3$
 \\
$G(1,7)$ & $1$ & $A^7B^7(A-B)^7(A^3-8A^2B+5AB^2+B^3)$\\
   & $7$ & $AB(A-B)(A^3-8A^2B+5AB^2+B^3)^7$\\ 
$G(1,8)$ & $1$ & $A^{8} B^{2} (A - B)^{8}(2 A - B)^{4}  (8 A^{2} - 8 A B + B^{2})
$\\
   & $2$ & $A^{4} B^{4} (A - B)^{4}  (2 A - B)^{8}  (8 A^{2} - 8 A B + B^{2})^{2}$\\
$G(1,9)$ & $1$ & $A^{9} B^{9}    (A - B)^{9}  (A^{2} - A B + B^{2})^{3}  (A^{3} - 6 A^{2} B + 3 A B^{2} + B^{3})
$\\
			  & $3$ & $A^{3}B^{3} (A - B)^{3} (A^{2} - A B + B^{2})^{9} (A^{3} - 6 A^{2} B + 3 A B^{2} + B^{3})^{3}
$
  \\ 
$G(1,10)$ & $1$ & $-A^{10} B^{5}  (A - B)^{10} (2 A - B)^{5} (A^{2} - 3 A B + B^{2})^{2} (4 A^{2} - 2 A B - B^{2}) 
$\\
			  & $2$ & $A^{5}B^{10}  (A - B)^{5}    (2 A - B)^{10} (A^{2} - 3 A B + B^{2}) (4 A^{2} - 2 A B - B^{2})^{2}$\\ 
			  & $5$ & $A^{2}B   (A - B)^{2}   (2 A - B)  (A^{2} - 3 A B + B^{2})^{10} (4 A^{2} - 2 A B - B^{2})^{5} 
$
   \\ 
$G(1,12)$ & $1$ & $A^{12}B^{2}   (A - B)^{12} (2 A - B)^{6}  (2 A^{2} - 2 A B + B^{2})^{3}(3 A^{2} - 3 A B + B^{2})^{4}(6 A^{2} - 6 A B + B^{2})    $\\
			  & $2$ & $ A^{6}B^{4}  (A - B)^{6}   (2 A - B)^{12}  (2 A^{2} - 2 A B + B^{2})^{6}(3 A^{2} - 3 A B + B^{2})^{2}  (6 A^{2} - 6 A B + B^{2})^{2}  $\\ 
			  & $3$ & $A^{4}B^{6}    (A - B)^{4} (2 A - B)^{2}   (2 A^{2} - 2 A B + B^{2})  (3 A^{2} - 3 A B + B^{2})^{12}(6 A^{2} - 6 A B + B^{2})^{3}  
$
  \\
$G(2,2)$ & $1$ & $A^2B^2(A-B)^2$\\
 & $2$ & $AB^4(A-B)$
  \\
$G(2,4)$ & $1$ & $2^{12} A^{2} B^{2}   (4 A - B)^{4}  (4 A + B)^{4}$\\
   & $2$ & $2^{24}   A^{4} B^{4}  (4 A - B)^{2}  (4 A + B)^{2} 
$
  \\
$G(2,6)$ & $1$ & $2^{12}   A^{6} B^{2} (A - B)^{6}  (A + B)^{6} (3 A - B)^{2}  (3 A + B)^{2}  
$\\
 & $2$ & $-2^{12}  A^{3} B (A - B)^{12} (A + B)^{3}  (3 A - B)  (3 A + B)^{4}  
$\\
 & $3$  & $2^{12} A^{2}  B^{6}  (A - B)^{2}  (A + B)^{2}  (3A -  B)^{6}  (3A +  B)^{6}$
  \\
$G(2,8)$ & $1$ & $2^{12} A^{8} B^{8}    (A - B)^{8}  (A + B)^{8}  (A^{2} - 2 A B - B^{2})^{2}  (A^{2} + 2 A B - B^{2})^{2}  (A^{2} + B^{2})^{4}
$\\
 & $2$ & $2^{12} A^{4} B^{4}    (A - B)^{4}  (A + B)^{4}  (A^{2} - 2 A B - B^{2})^{4}  (A^{2} + 2 A B - B^{2})^{4}  (A^{2} + B^{2})^{8} $
  \\
$G(3,3)$ & $1$ & $6^{12} B^{3} (A - 3 B)^{3}  (A^{2} + 3 A B + 9 B^{2})^{3}$\\
 & 3 & $6^{12} B (A - 3 B)^{9}  (A^{2} + 3 A B + 9 B^{2})
$
  \\
$G(3,6)$ & $1$ & $3^6 A^{6} B^{6} (3 A + B)^{3} (9 A + B)^{3} (27 A^{2} + B^{2})^{3} (27 A^{2} + 9 A B + B^{2})^{6}  
$\\
 & $2$ & $-3^{3}   A^{3} B^{3} (3 A + B)^{6}  (9 A + B)^{6}  (27 A^{2} + B^{2})^{6} (27 A^{2} + 9 A B + B^{2})^{3} 
$\\
 & $3$ & $3^{2}   A^{2} B^{18} (3 A + B)  A^{2}  (9 A + B)^{9}    (27 A^{2} + B^{2})  (27 A^{2} + 9 A B + B^{2})^{2}
$
  \\
$G(4,4)$ & $1$ & $2^8 A^{4} B^{4}    (A - B)^{4}  (2 A - B)^{4}  (2 A^{2} - 2 A B + B^{2})^{4}
$\\
 & $2$ & $2^{10}  A^{2}  B^{8}  (A - B)^{2}   (2 A - B)^{8}  (2 A^{2} - 2 A B + B^{2})^{2}
$
  \\
$G(5,5)$ & $1$ & $6^{12} A^{5}  B^{5}   (A^{2} - A B - B^{2})^{5}  (A^{4} - 2 A^{3} B + 4 A^{2} B^{2} - 3 A B^{3} + B^{4})^{5}  (A^{4} + 3 A^{3} B + 4 A^{2} B^{2} + 2 A B^{3} + B^{4})^{5}
$\\
& $5$ & $6^{12} A  B   (A^{2} - A B - B^{2})^{25}  (A^{4} - 2 A^{3} B + 4 A^{2} B^{2} - 3 A B^{3} + B^{4})  (A^{4} + 3 A^{3} B + 4 A^{2} B^{2} + 2 A B^{3} + B^{4})
$\\
\bottomrule[1.5pt]
\end{longtblr}


\begin{longtblr}[
   caption = {Constants for asymptotics for torsion families.}, 
   label = {tab:elliptic_count}
]{
   width = \textwidth,  
   rowhead = 1,
   row{1} = {white},
   colspec = {X[c] c X[c] X[c] c}
}
\toprule[1.5pt]
 $G$  & $(w_0,w_1)$ & $\alpha(T)=\alpha(G)$ & $\beta(T)=\beta(G)$\\
\midrule
 $G(1,1)$ & $(4,6)$ & $5/6$ & $\frac{5}{6}-\frac{1}{3\deg_{\Q}(K)}$\\
 $G(1,2)$ & $(2,4)$ & $1/2$ & $\frac{1}{2}-\frac{1}{6\deg_{\Q}(K)}$\\
 $G(1,3)$ & $(1,3)$ & $1/3$ & $\frac{1}{3}-\frac{1}{12\deg_{\Q}(K)}$\\
 $G(1,4)$ & $(1,2)$ & $1/4$ & $\frac{1}{4}-\frac{1}{12\deg_{\Q}(K)}$\\
 $G(1,5)$ & $(1,1)$ & $1/6$ & $\frac{1}{6}-\frac{1}{12\deg_{\Q}(K)}$\\
 $G(1,6)$ & $(1,1)$ & $1/6$ & $\frac{1}{6}-\frac{1}{12\deg_{\Q}(K)}$\\
 $G(1,7)$ & $(1,1)$ & $1/12$ & $\frac{1}{12}-\frac{1}{24\deg_{\Q}(K)}$\\
 $G(1,8)$ & $(1,1)$ & $1/12$ & $\frac{1}{12}-\frac{1}{24\deg_{\Q}(K)}$\\
$G(1,9)$ & $(1,1)$ & $1/18$ & $\frac{1}{18}-\frac{1}{36\deg_{\Q}(K)}$\\
 $G(1,10)$ & $(1,1)$ & $1/18$ & $\frac{1}{18}-\frac{1}{36\deg_{\Q}(K)}$\\
 $G(1,12)$ & $(1,1)$ &  $1/24$ & $\frac{1}{24}-\frac{1}{48\deg_{\Q}(K)}$\\
 $G(2,2)$  & $(2,2)$ & $1/3$ & $\frac{1}{3}-\frac{1}{6\deg_{\Q}(K)}$\\
 $G(2,4)$ & $(1,1)$ & $1/6$ & $\frac{1}{6}-\frac{1}{12\deg_{\Q}(K)}$\\
  $G(2,4)$  &  $(1,1)$ & $1/12$ & $\frac{1}{12}-\frac{1}{24\deg_{\Q}(K)}$\\
  $G(2,8)$ & $(1,1)$ & $1/24$ & $\frac{1}{24}-\frac{1}{48\deg_{\Q}(K)}$\\
  $G(3,3)$ & $(1,1)$ & $1/6$ & $\frac{1}{6}-\frac{1}{12\deg_{\Q}(K)}$\\
  $G(3,6)$ & $(1,1)$ & $1/18$ & $\frac{1}{18}-\frac{1}{36\deg_{\Q}(K)}$\\
  $G(4,4)$ & $(1,1)$ & $1/12$ & $\frac{1}{12}-\frac{1}{24\deg_{\Q}(K)}$\\
  $G(5,5)$ & $(1,1)$ & $1/30$ & $\frac{1}{30}-\frac{1}{60\deg_{\Q}(K)}$\\
\bottomrule[1.5pt]
\end{longtblr}


\begin{longtblr}[ caption = {Local densities for elliptic curves with prescribed level structure at finite places not above $2$ or $3$.}, label = {tab:LocalConditions},]{
  width = \linewidth,  rowhead = 2,
  row{4, 7-8, 13-14, 19-22, 27-30, 35, 38-39, 44-45, 48-49} = {gray!20}, row{1} = {white}, 
colspec={X[c]cX[c]X[c]X[c]X[c]}}
\toprule[1.5pt]
  Level Structure & Conditions on & good & additive & mult. & split mult.\\
  $G$ & $N_{K/\Q}(\p_v)=q$ & $\kappa_{good}(G,v)$ & $\kappa_{add}(G,v)$ & $\kappa_{mult}(G,v)$ & $\kappa_{split}(G,v)$\\
\midrule
$G(1,1)$ &  &$\frac{q-1}{q}$ & $\frac{q^8-1}{q^{10}}$ & $\frac{q-1}{q^2}$ & $\frac{q-1}{2q^2}$ \\
$G(1,2)$ & & $\frac{q(q-1)}{q^2}$ & $\frac{q^4-1}{q^6}$ & $\frac{2q-2}{q^2}$ & $\frac{q-1}{q^2}$\\
%
%
$G(1,3)$ & $q\equiv 1\pmod{3}$ & $\frac{(q-1)^2}{q^2}$ & $\frac{q^2-1}{q^4}$ & $\frac{2q-2}{q^2}$ & $\frac{2q-2}{q^2}$\\
 & $q\equiv 2\pmod{3}$ & $\frac{(q-1)^2}{q^2}$ & $\frac{q^2-1}{q^4}$ & $\frac{2q-2}{q^2}$ & $\frac{q-1}{q^2}$\\
$G(1,4)$ & $q\equiv 1 \pmod{4}$ & $\frac{(q-1)(q-2)}{q^2}$ & $\frac{q-1}{q^3}$ & $\frac{3q-3}{q^2}$ & $\frac{5q-5}{2q^2}$ \\
& $q\equiv 3 \pmod{4}$ & $\frac{(q-1)(q-2)}{q^2}$ & $\frac{q-1}{q^3}$ & $\frac{3q-3}{q^2}$ & $\frac{3q-3}{2q^2}$ \\
$G(1,5)$ & $q\equiv  0 \pmod{5}$ & $\frac{(q-1)(q-2)}{q^2}$ & $\frac{q-1}{q^2}$ & $\frac{2q-2}{q^2}$  & $\frac{2q-2}{q^2}$ \\
& $q\equiv  1 \pmod{5}$ & $\frac{(q-1)(q-3)}{q^2}$ & $0$ & $\frac{4q-4}{q^2}$ & $\frac{4q-4}{q^2}$\\
& $q\equiv  4 \pmod{5}$ & $\frac{(q-1)(q-3)}{q^2}$ & $0$ & $\frac{4q-4}{q^2}$ & $\frac{2q-2}{q^2}$\\
& $q\equiv  \pm 2 \pmod{5}$ & $\frac{(q-1)^2}{q^2}$ & $0$ & $\frac{2q-2}{q^2}$ & $\frac{2q-2}{q^2}$ \\
%
%
$G(1,6)$ &  $q\equiv 1\pmod{6}$ & $\frac{(q-1)(q-3)}{q^2}$ & $0$ & $\frac{4q-4}{q^2}$ & $\frac{4q-4}{q^2}$\\
 & $q\equiv 5 \pmod{6}$ & $\frac{(q-1)(q-3)}{q^2}$ & $0$ & $\frac{4q-4}{q^2}$ & $\frac{2q-2}{q^2}$\\
%
%
%
$G(1,7)$ & $q\equiv 0 \pmod{7}$ & $\frac{(q-1)(q-3)}{q^2}$ & $\frac{q-1}{q^2}$ & $\frac{3q-3}{q^2}$ & $\frac{3q-3}{q^2}$  \\
%
%
&$q\equiv 1 \pmod{7}$  & $\frac{(q-1)(q-5)}{q^2}$ & $0$ & $\frac{6q-6}{q^2}$ & $\frac{6q-6}{q^2}$ \\
&$q\equiv 2,3,4,5 \pmod{7}$  & $\frac{(q-1)(q-2)}{q^2}$ & $0$ & $\frac{3q-3}{q^2}$ & $\frac{3q-3}{q^2}$ \\
&$q\equiv 6 \pmod{7}$  & $\frac{(q-1)(q-5)}{q^2}$ & $0$ & $\frac{6q-6}{q^2}$ & $\frac{3q-3}{q^2}$ \\
%
%
$G(1,8)$ & $q\equiv 1 \pmod{8}$  & $\frac{(q-1)(q-5)}{q^2}$ & $0$ & $\frac{6q-6}{q^2}$ & $\frac{6q-6}{q^2}$ \\
& $q\equiv - 1 \pmod{8}$  & $\frac{(q-1)(q-5)}{q^2}$ & $0$ & $\frac{6q-6}{q^2}$ & $\frac{3q-3}{q^2}$ \\
 & $q\equiv 3 \pmod{8}$  & $\frac{(q-1)(q-3
)}{q^2}$ & $0$ & $\frac{4q-4}{q^2}$ & $\frac{3q-3}{q^2}$\\
 & $q\equiv - 3 \pmod{8}$  & $\frac{(q-1)(q-3
)}{q^2}$ & $0$ & $\frac{4q-4}{q^2}$ & $\frac{4q-4}{q^2}$ \\
%
%
%
%
$G(1,9)$ & $q\equiv 1 \pmod{9}$  & $\frac{(q-1)(q-7)}{q^2}$ & $0$ & $\frac{8q-8}{q^2}$ & $\frac{8q-8}{q^2}$ \\
&$q\equiv 2, 5 \pmod{9}$  & $\frac{(q-1)(q-2)}{q^2}$ & $0$ & $\frac{3q-3}{q^2}$ & $\frac{3q-3}{q^2}$ \\
&$q\equiv  4, 7 \pmod{9}$  & $\frac{(q-1)(q-4)}{q^2}$ & $0$ & $\frac{5q-5}{q^2}$ & $\frac{5q-5}{q^2}$ \\
&$q\equiv 8 \pmod{9}$  & $\frac{(q-1)(q-5)}{q^2}$ & $0$ & $\frac{6q-6}{q^2}$ & $\frac{3q-3}{q^2}$  \\
%
%
%
%
$G(1,10)$ & $q\equiv 5 \pmod{10}$  & $\frac{(q-1)(q-4)}{q^2}$ & $\frac{q-1}{q^2}$ & $\frac{4q-4}{q^2}$ & $\frac{4q-4}{q^2}$ \\
&$q\equiv  1 \pmod{10}$  & $\frac{(q-1)(q-7)}{q^2}$ & $0$ & $\frac{8q-8}{q^2}$ & $\frac{8q-8}{q^2}$ \\
&$q\equiv 9 \pmod{10}$  & $\frac{(q-1)(q-7)}{q^2}$ & $0$ & $\frac{8q-8}{q^2}$ & $\frac{4q-4}{q^2}$ \\
&$q\equiv \pm 3 \pmod{10}$  & $\frac{(q-1)(q-3)}{q^2}$ & $0$ & $\frac{4q-4}{q^2}$ & $\frac{4q-4}{q^2}$ \\
%
%
%
%
$G(1,12)$ & $q\equiv 1 \pmod{12}$  & $\frac{(q-1)(q-9)}{q^2}$ & $0$ & $\frac{10q-10}{q^2}$ & $\frac{10q-10}{q^2}$ \\
& $q\equiv 5 \pmod{12}$  & $\frac{(q-1)(q-5)}{q^2}$ & $0$ & $\frac{6q-6}{q^2}$  & $\frac{5q-5}{q^2}$\\
& $q\equiv 7 \pmod{12}$  & $\frac{(q-1)(q-5)}{q^2}$ & $0$ & $\frac{6q-6}{q^2}$  & $\frac{6q-6}{q^2}$ \\
& $q\equiv 11 \pmod{12}$  & $\frac{(q-1)(q-5)}{q^2}$ & $0$ & $\frac{6q-6}{q^2}$  & $\frac{3q-3}{q^2}$\\
$G(2,2)$ & & $\frac{(q-1)(q-2)}{q^2}$ & $\frac{q^2-1}{q^4}$ & $\frac{3q-3}{q^2}$ & $\frac{3q-3}{2q^2}$ \\
$G(2,4)$ & $q\equiv 1\pmod{4}$ & $\frac{(q-1)(q-3)}{q^2}$ & $0$ & $\frac{4q-4}{q^2}$ & $\frac{4q-4}{q^2}$ \\
& $q\equiv 3\pmod{4}$ & $\frac{(q-1)(q-3)}{q^2}$ & $0$ & $\frac{4q-4}{q^2}$ & $\frac{2q-2}{q^2}$ \\
$G(2,6)$ & $q\equiv 1\pmod{3}$ & $\frac{(q-1)(q-5)}{q^2}$ & $0$ & $\frac{6q-6}{q^2}$ & $\frac{6q-6}{q^2}$ \\
& $q\equiv 2\pmod{3}$ & $\frac{(q-1)(q-3)}{q^2}$ & $0$ & $\frac{6q-6}{q^2}$ & $\frac{3q-3}{q^2}$ \\
$G(2,8)$ & $q\equiv 1\pmod{8}$ & $\frac{(q-1)(q-9)}{q^2}$ & $0$ & $\frac{10q-10}{q^2}$ & $\frac{10q-10}{q^2}$ \\
& $q\equiv 3\pmod{8}$ & $\frac{(q-1)(q-3)}{q^2}$ & $0$ & $\frac{4q-4}{q^2}$ & $\frac{4q-4}{q^2}$ \\
& $q\equiv 5\pmod{8}$ & $\frac{(q-1)(q-5)}{q^2}$ & $0$ & $\frac{6q-6}{q^2}$ & $\frac{6q-6}{q^2}$ \\
& $q\equiv 7\pmod{8}$ & $\frac{(q-1)(q-7)}{q^2}$ & $0$ & $\frac{8q-8}{q^2}$ & $\frac{4q-4}{q^2}$ \\
$G(3,3)$ &  $q\equiv 1\pmod{3}$ & $\frac{(q-1)(q-3)}{q^2}$ & $0$ & $\frac{4q-4}{q^2}$ & $\frac{4q-4}{q^2}$ \\
 &  $q\equiv 2\pmod{3}$ & $\frac{(q-1)^2}{q^2}$ & $0$ & $\frac{2q-2}{q^2}$ & $\frac{q-1}{q^2}$\\
$G(3,6)$ &  $q\equiv 1\pmod{6}$ & $\frac{(q-1)(q-7)}{q^2}$ & $0$ & $\frac{8q-8}{q^2}$ & $\frac{8q-8}{q^2}$ \\
&  $q\equiv 5\pmod{6}$ & $\frac{(q-1)(q-3)}{q^2}$ & $0$ & $\frac{4q-4}{q^2}$ & $\frac{2q-2}{q^2}$ \\
$G(4,4)$ &  $q\equiv 1\pmod{4}$ & $\frac{(q-1)(q-5)}{q^2}$ & $0$ & $\frac{6q-6}{q^2}$ & $\frac{6q-6}{q^2}$ \\
& $q\equiv 3\pmod{4}$ & $\frac{(q-1)(q-3)}{q^2}$ & $0$ & $\frac{4q-4}{q^2}$ & $\frac{2q-2}{q^2}$ \\
$G(5,5)$ &  $q\equiv 1\pmod{5}$ & $\frac{(q-1)(q-11)}{q^2}$ & $0$ & $\frac{12q-12}{q^2}$ & $\frac{12q-12}{q^2}$ \\
&  $q\equiv 2,3 \pmod{5}$ & $\frac{(q-1)^2}{q^2}$ & $0$ & $\frac{2q-2}{q^2}$ & $\frac{2q-2}{q^2}$ \\
%
%
&  $q\equiv 4\pmod{5}$ & $\frac{(q-1)(q-3)}{q^2}$ & $0$ & $\frac{4q-4}{q^2}$ & $\frac{2q-2}{q^2}$ \\
\bottomrule[1.5pt]
\end{longtblr}


\begin{longtblr}[ caption = {Local densities for $n$-reduction with respect to isogenies $\phi$ of prime degree at places not above $2$, $3$, or $\deg(\phi)$.}, label = {tab:isogeny_local_densities}]{
  width = \linewidth, rowhead = 2,
  row{3,6,10-12,16-19, 24-28, 36, 39-42, 47-48, 53-54} = {gray!20}, row{1} = {white},  
colspec={X[c]ccc},}
\toprule[1.5pt]
  Level Structure & Isogeny & Conditions on & $n$-reduction at $\phi$  \\
$G$ & $\Phi$ & $N_{K/\Q}(\p_v)=q$ & $\left(\kappa_{-p}(G,\phi_\Phi, v),\kappa_{0}(G,\phi_\Phi, v),\kappa_{p}(G,\phi_\Phi, v)\right)$ \\
\midrule
%
%
$G(1,2)$ & $C_2$ &  
& $\left(\frac{q-1}{q^2},0,\frac{q-1}{q^2}\right)$ \\
%
%
$G(1,3)$ & $C_3$ & $q\equiv 1\pmod{3}$ & $\left(\frac{q-1}{q^2},0,\frac{q-1}{q^2}\right)$  \\
 &  & $q\equiv 2\pmod{3}$ & $\left(\frac{q-1}{q^2},\frac{q-1}{q^2},0\right)$  \\
%
%
$G(1,4)$ & $C_2$ &  & $\left(\frac{q-1}{q^2},\frac{q-1}{2q^2},\frac{3(q-1)}{2q^2}\right)$ \\
%
%
$G(1,5)$ & $C_5$ & $q\equiv  1 \pmod{5}$ & $\left(\frac{2q-2}{q^2},0,\frac{2q-2}{q^2}\right)$  \\
&  & $q\equiv 4 \pmod{5}$ & $\left(\frac{2q-2}{q^2},\frac{2q-2}{q^2},0\right)$  \\
&  & $q\equiv 2,3 \pmod{5}$ & $\left(\frac{2q-2}{q^2},0,0\right)$  \\
%
%
$G(1,6)$ &  $C_2$ &  & $\left(\frac{2q-2}{q^2},0,\frac{2q-2}{q^2}\right)$\\
 &  $C_3$ & $q\equiv 1\pmod{6}$ & $\left(\frac{2q-2}{q^2},0,\frac{2q-2}{q^2}\right)$\\
 &  & $q\equiv 5 \pmod{6}$ & $\left(\frac{2q-2}{q^2},\frac{2q-2}{q^2},0\right)$\\
%
%
%
$G(1,7)$ & $C_7$ &  $q\equiv 1 \pmod{7}$  & $\left(\frac{3q-3}{q^2},0,\frac{3q-3}{q^2}\right)$ \\
& & $q\equiv 2,3,4,5 \pmod{7}$  & $\left(\frac{3q-3}{q^2},0, 0\right)$ \\
 &  &  $q\equiv 6 \pmod{7}$  & $\left(\frac{3q-3}{q^2},\frac{3q-3}{q^2},0 \right)$ \\
%
%
%
%
$G(1,8)$ & $C_2$ & $q\equiv 1 \pmod{8}$  
& $\left(\frac{2q-2}{q^2}, 0, \frac{4q-4}{q^2}\right)$ \\
& & $q\equiv 3 \pmod{8}$  
& $\left(\frac{2q-2}{q^2}, \frac{q-1}{q^2}, \frac{q-1}{q^2}\right)$ \\
& & $q\equiv 5 \pmod{8}$  
& $\left(\frac{2q-2}{q^2}, 0, \frac{2q-2}{q^2}\right)$ \\
& & $q\equiv 7 \pmod{8}$  
& $\left(\frac{2q-2}{q^2}, \frac{q-1}{q^2}, \frac{3q-3}{q^2}\right)$ \\
%
%
%
%
$G(1,9)$ & $C_3$ & $q\equiv 1 \pmod{9}$  & $\left(\frac{3q-3}{q^2},0,\frac{5q-5}{q^2}\right)$ \\
& & $q\equiv 2, 5 \pmod{9}$  & $\left(\frac{3q-3}{q^2},0,0\right)$ \\
& & $q\equiv  4, 7 \pmod{9}$  &  $\left(\frac{3q-3}{q^2},0,\frac{2q-2}{q^2}\right)$ \\
& & $q\equiv 8 \pmod{9}$  &  $\left(\frac{3q-3}{q^2},\frac{3q-3}{q^2},0\right)$  \\
%
%
%
%
$G(1,10)$ & $C_2$ & $q\equiv 1, 9 \pmod{10}$  &  $\left(\frac{4q-4}{q^2},0,\frac{4q-4}{q^2}\right)$ \\
&  & $q\equiv 3, 7 \pmod{10}$  &  $\left(\frac{2q-2}{q^2},0,\frac{2q-2}{q^2}\right)$ \\
& $C_5$ & $q\equiv 1 \pmod{10}$  & $\left(\frac{4q-4}{q^2},0,\frac{4q-4}{q^2}\right)$  \\
& & $q\equiv 9 \pmod{10}$  & $\left(\frac{4q-4}{q^2},\frac{4q-4}{q^2},0\right)$  \\
& & $q\equiv 3, 7 \pmod{10}$  & $\left(\frac{4q-4}{q^2},0,0\right)$  \\
%
%
%
%
$G(1,12)$ & $C_2$ & $q\equiv 1 \pmod{12}$  & $\left(\frac{4q-4}{q^2},0,\frac{6q-6}{q^2}\right)$ \\
& & $q\equiv 5, 11 \pmod{12}$  & $\left(\frac{2q-2}{q^2},\frac{q-1}{2},\frac{3q-3}{q^2}\right)$ \\
& & $q\equiv 7 \pmod{12}$  & $\left(\frac{4q-4}{q^2},0,\frac{2q-2}{q^2}\right)$ \\
 & $C_3$ & $q\equiv 1 \pmod{12}$  & $\left(\frac{5q-5}{q^2},0,\frac{5q-5}{q^2}\right)$ \\
& & $q\equiv 5 \pmod{12}$  & $\left(\frac{5q-5}{q^2},\frac{q-1}{q^2},0\right)$ \\
& & $q\equiv 7 \pmod{12}$  & $\left(\frac{3q-3}{q^2},0,\frac{3q-3}{q^2}\right)$ \\
& & $q\equiv 11 \pmod{12}$  & $\left(\frac{3q-3}{q^2},\frac{3q-3}{q^2},0\right)$ \\
%
%
%
%
$G(2,2)$  &  $C_2$ &  & $\left(\frac{2q-2}{q^2},\frac{q-1}{2q^2},\frac{q-1}{2q^2}\right)$ \\
%
%
%
%
$G(2,4)$  &  $C_2$ & $q\equiv 1\pmod{4}$ & $\left(\frac{2q-2}{q^2},0,\frac{2q-2}{q^2}\right)$ \\
 &  & $q\equiv 3\pmod{4}$ & $\left(\frac{2q-2}{q^2},\frac{2q-2}{q^2},0\right)$ \\
%
%
%
%
$G(2,6)$  &  $C_2$ & $q\equiv 1\pmod{6}$ & $\left(\frac{4q-4}{q^2},0,\frac{2q-2}{q^2}\right)$ \\
 &  & $q\equiv 5\pmod{6}$ & $\left(\frac{4q-4}{q^2},\frac{q-1}{q^2},\frac{q-1}{q^2}\right)$ \\
 & $C_3$ & $q\equiv 1 \pmod{6}$  & $\left(\frac{3q-3}{q^2},\frac{q-1}{q^2},\frac{2q-2}{q^2}\right)$ \\
& & $q\equiv 5 \pmod{6}$  & $\left(\frac{3q-3}{q^2},\frac{3q-3}{q^2},0\right)$ \\
%
%
%
%
$G(2,8)$  &  $C_2$ & $q\equiv 1\pmod{8}$ & $\left(\frac{4q-4}{q^2},0,\frac{6q-6}{q^2}\right)$ \\
 &  & $q\equiv 3\pmod{8}$ & $\left(\frac{4q-4}{q^2},0,0\right)$ \\
 &  & $q\equiv 5\pmod{8}$ & $\left(\frac{4q-4}{q^2},0,\frac{2q-2}{q^2}\right)$ \\
 &  & $q\equiv 7\pmod{8}$ & $\left(\frac{4q-4}{q^2},\frac{4q-4}{q^2},0\right)$ \\
%
%
%
%
$G(3,3)$ &  $C_3$ & $q\equiv 1\pmod{3}$ & $\left(\frac{3q-3}{q^2}, 0, \frac{q-1}{q^2}\right)$ \\
 & &  $q\equiv 2\pmod{3}$ & $\left(\frac{q-1}{q^2}, \frac{q-1}{q^2}, 0\right)$ \\
%
%
%
%
$G(3,6)$ & $C_2$ & $q\equiv 1\pmod{6}$ & $\left(\frac{4q-4}{q^2}, 0, \frac{4q-4}{q^2}\right)$ \\ 
& &  $q\equiv 5\pmod{6}$ & $\left(\frac{2q-2}{q^2}, 0, \frac{2q-2}{q^2}\right)$ \\ 
& $C_3$ & $q\equiv 1\pmod{6}$ & $\left(\frac{6q-6}{q^2}, 0, \frac{2q-2}{q^2}\right)$ \\ 
& &  $q\equiv 5\pmod{6}$ & $\left(\frac{2q-2}{q^2}, \frac{2q-2}{q^2}, 0\right)$ \\ 
%
%
%
%
$G(4,4)$ & $C_2$ &  $q\equiv 1\pmod{4}$ & $\left(\frac{4q-4}{q^2}, 0, \frac{2q-2}{q^2}\right)$ \\
& & $q\equiv 3\pmod{4}$ & $\left(\frac{2q-2}{q^2}, \frac{2q-2}{q^2}, 0\right)$ \\ 
%
%
%
%
$G(5,5)$ & $C_5$ &  $q\equiv 1\pmod{5}$ & $\left(\frac{10q-10}{q^2}, 0, \frac{2q-2}{q^2}\right)$ \\ 
& &  $q\equiv 2, 3\pmod{5}$ & $\left(\frac{2q-2}{q^2}, 0, 0\right)$ \\ 
%
%
& & $q\equiv 4\pmod{5}$ & $\left(\frac{2q-2}{q^2}, \frac{2q-2}{q^2}, 0\right)$ \\ 
\bottomrule[1.5pt]
\end{longtblr}


\begin{longtblr}[ caption = {Mean and variance for isogenies of prime degree I.},label = {tab:mean_variance_I}]{
  width = \linewidth, rowhead = 1,
  row{2, 5-6, 10-13, 18-22, 27-32, 43, 46-49, 55-56, 61-62} = {gray!20}, row{1} = {white}, 
colspec={ccccccc}, 
}
\toprule[1.5pt]
$G$ & $m$ & $\Phi$ &  $\cG_m(K)$ &  $c_{-}(G,\Phi,\cG_m(K))$ & $ c_{+}(G,\Phi,\cG_m(K))$ \\
\midrule
$G(1,2)$ & $2$ & $C_2$ & $\{1\}$ & $1$ & $1$ \\
%
%
$G(1,3)$ &  $3$ & $C_3$ & $\{1\}$ & $1$ &  $1$ \\
 &  &  & $(\Z/3\Z)^\times$ & $1$ &  $1/2$ \\
%
%
%
%
$G(1,4)$ &  $4$ & $C_2$ & $\{1\}$ & $1$ &  $3/2$ \\
 &  &  & $(\Z/4\Z)^\times$ & $1$ &  $3/2$  \\
%
%
%
%
$G(1,5)$ &  $5$ & $C_5$ & $\{1\}$ & $2$ &  $2$ \\
 &  &  & $\{1,4\}$  &  $2$ & $1$ \\
 &  &  & $(\Z/5\Z)^\times$ & $2$ &  $1/2$ \\
 %
%
%
%
    $G(1,6)$ &  $6$ & $C_2$ & $\{1\}$ & $2$ & $2$ \\
    & & & $(\Z/6\Z)^\times$ & $2$ & $2$ \\
     &   & $C_3$ & $\{1\}$ &  $2$ & $2$ \\
     &   &   & $(\Z/6\Z)^\times$ &   $2$ & $1$ \\
%
%
%
%
    $G(1,7)$ & $7$ &  $C_7$  & $\{1\}$ & $3$ &  $3$\\
     &  &    & $\{1,6\}$ & $3$ & $3/2$ \\
    &  &    & $\{1, 2, 4\}$ & $3$ & $1$\\
    &  &    & $(\Z/7\Z)^\times$ & $3$ & $1/2$ \\
    %
    %
    %
    %
    %
    %
    %
    %
    %
%
%
%
%
$G(1,8)$ &  $8$ &  $C_2$  & $\{1\}$ &  $2$ & $4$ \\
     &  &  & $\{1,3\}$ &  $2$ & $5/2$ \\
    &  &  & $\{1,5\}$ &  $2$ & $3$ \\
     &  &  & $\{1,7\}$ &  $2$ & $7/2$ \\
     &  &  & $(\Z/8\Z)^\times$ &  $2$ & $5/2$ \\
%
%
%
%
$G(1,9)$ &  $9$ & $C_3$ & $\{1\}$ & $3$ & $5$ \\
 &  &  & $\{1,8\}$ & $3$ & $5/2$ \\
 &  &  & $\{1,4,7\}$ & $3$ &  $3$\\
 &  &  & $(\Z/9\Z)^\times$ & $3$ & $3/2$ \\
%
%
%
%
$G(1,10)$ &  $10$ & $C_2$ & $\{1\}$ & $4$  & $4$ \\
 &  &  & $\{1,9\}$ & $4$ & $4$\\
 &  &  & $(\Z/10\Z)^\times$ &  $3$ & $3$ \\
&   & $C_5$ & $\{1\}$ & $4$ & $4$\\
 &  &  & $\{1,9\}$ & $4$ & $2$\\
 &  &  & $(\Z/10\Z)^\times$ & $4$ & $1$ \\
%
%
%
%
$G(1,12)$ & $12$ & $C_2$ & $\{1\}$ & $4$ & $6$\\
 &  &  & $\{1,5\}$ & $3$ & $9/2$\\
 &  &  & $\{1,7\}$ & $4$ & $4$\\
 &  &  & $\{1,11\}$ & $3$ & $9/2$\\
 &  &  & $(\Z/12\Z)^\times$ & $3$ & $7/2$ \\
& & $C_3$  & $\{1\}$ & $5$ & $5$\\
 &  &  & $\{1,5\}$ & $5$ & $5/2$\\
 &  &  & $\{1,7\}$ & $4$ & $4$\\
 &  &  & $\{1,11\}$ & $4$ & $5/2$ \\
 &  &  & $(\Z/12\Z)^\times$ & $4$ & $2$\\
%
%
%
%
$G(2,2)$ & $2$ & $C_2$ & $\{1\}$ & $2$ & $1/2$ \\
%
%
%
%
    $G(2,4)$ & $4$ &  $C_2$ & $\{1\}$ & $2$ & $2$ \\
    &   & & $(\Z/4\Z)^\times$ & $2$ & $1$\\
%
%
%
%
    $G(2,6)$ & $6$ & $C_2$ &   $\{1\}$ & $4$ & $2$ \\
    &  &  & $(\Z/6\Z)^\times$ & $4$ & $3/2$ \\
     &  & $C_3$ &  $\{1\}$ & $3$ & $2$ \\
    &  &  & $(\Z/6\Z)^\times$ & $3$ & $1$ \\ 
%
%
%
%
    $G(2,8)$ & $8$ & $C_2$ &   $\{1\}$ & $4$ & $6$\\
    &  &  & $\{1, 3\}$ & $4$ & $3$\\
    &  &  & $\{1, 5\}$ & $4$ & $4$\\
    &  &  & $\{1, 7\}$ & $4$ & $3$\\
    &  &  & $(\Z/8\Z)^\times$ & $4$ & $2$\\
%
%
%
%
    $G(3,3)$ & $3$ & $C_3$ & $\{1\}$ & $3$ & $1$\\
    &  &  & $(\Z/3\Z)^\times$ & $2$ & $1/2$\\
%
%
%
%
    $G(3,6)$ & $6$ & $C_2$ & $\{1\}$ & $4$ & $4$ \\
    &  & & $(\Z/6\Z)^\times$ & $3$ & $3$ \\
    &  & $C_3$ & $\{1\}$ & $6$ & $2$ \\
    &  & & $(\Z/6\Z)^\times$ & $4$ & $1$ \\
%
%
%
%
    $G(4,4)$ & $2$ & $C_2$ & $\{1\}$ &  $4$ & $2$ \\
    & & & $(\Z/4\Z)^\times$ & $3$ & $1$\\
%
%
%
%
    $G(5,5)$ & $5$ & $C_5$ & $\{1\}$ & $10$ & $2$ \\
    & &  & $\{1,4\}$ & $6$ & $1$ \\
     & &  & $(\Z/5\Z)^\times$ & $4$ & $1/2$ \\
\bottomrule[1.5pt]
\end{longtblr}


\begin{longtblr}[ caption = {Mean and variance for isogenies of prime degree II.}, label = {tab:mean_variance_II}]{
  width = \linewidth, rowhead = 1,
  row{2, 5-6, 10-13, 18-22, 27-32, 43, 46-49, 55-56, 61-62} = {gray!20}, row{1} = {white}, 
colspec={ccccccc}, 
}
\toprule[1.5pt]
$G$ & $m$ & $\Phi$ &  $\sG_m(K)$ &  $c_{\E}(G,\Phi,\sG_m(K))$ & $c_{\mathbb{V}}(G,\Phi,\sG_m(K))$ & $\theta(G,\Phi,\sG_m(K))$\\
\midrule
$G(1,2)$ & $2$ & $C_2$ & $\{1\}$ & $0$ &  $2$ & $1$ \\
%
%
$G(1,3)$ &  $3$ & $C_3$ & $\{1\}$ & $0$ &  $2$ & $1$ \\
 &  &  & $(\Z/3\Z)^\times$ & $-1/2$ &  $3/2$ & $1/4$ \\
$G(1,4)$ &  $4$ & $C_2$ & $\{1\}$ & $1/2$ &  $5/2$ & $7/4$\\
 &  &  & $(\Z/4\Z)^\times$ & $1/2$ &  $5/2$ & $7/4$ \\
%
%
%
%
$G(1,5)$ &  $5$ & $C_5$ & $\{1\}$ & $0$ &  $4$ & $2$ \\
 &  &  & $\{1,4\}$ & $-1$ &  $3$ & $1/2$ \\
 &  &  & $(\Z/5\Z)^\times$ & $-3/2$ &  $5/2$ & $-1/4$ \\
 %
%
%
%
    $G(1,6)$ &  $6$ & $C_2$ & $\{1\}$ & $0$ &  $4$ & $2$ \\
    & & & $(\Z/6\Z)^\times$ & $0$ & $4$ & $2$ \\
     &   & $C_3$ & $\{1\}$ & $0$ &  $4$ & $2$ \\
     &   &   & $(\Z/6\Z)^\times$ & $-1$ &  $3$ & $1/2$ \\
%
%
%
%
    $G(1,7)$ & $7$ &  $C_7$  & $\{1\}$ & $0$ & $6$ & $3$\\
     &  &    & $\{1,6\}$ & $-\frac{3}{2}$ & $\frac{9}{2}$ & $3/4$ \\
    &  &    & $\{1, 2, 4\}$ & $-2$ & $4$ & $0$\\
    &  &    & $(\Z/7\Z)^\times$ & $-\frac{5}{2}$ & $\frac{7}{2}$ & $-3/4$ \\
    %
    %
    %
    %
    %
    %
    %
    %
    %
%
%
%
%
$G(1,8)$ &  $8$ &  $C_2$  & $\{1\}$ &  $2$ & $6$ & $5$ \\
     &  &  & $\{1,3\}$ &  $1/2$ & $9/2$ & $11/4$ \\
    &  &  & $\{1,5\}$ &  $1$ & $5$ & $7/2$ \\
     &  &  & $\{1,7\}$ &  $3/2$ & $11/2$ & $17/4$ \\
     &  &  & $(\Z/8\Z)^\times$ &  $1/2$ & $9/2$ & $11/4$ \\
%
%
%
%
$G(1,9)$ &  $9$ & $C_3$ & $\{1\}$ & $2$ & $8$ & $6$ \\
 &  &  & $\{1,8\}$ & $-\frac{1}{2}$ & $\frac{11}{2}$ & $9/4$ \\
 &  &  & $\{1,4,7\}$ & $0$ & $6$ & $3$\\
 &  &  & $(\Z/9\Z)^\times$ & $-\frac{3}{2}$ & $\frac{9}{2}$ & $3/4$ \\
%
%
%
%
$G(1,10)$ &  $10$ & $C_2$ & $\{1\}$ & $0$ & $8$ & $4$ \\
 &  &  & $\{1,9\}$ & $0$ & $8$ & $4$\\
 &  &  & $(\Z/10\Z)^\times$ &  $0$ & $6$ & $3$ \\
&   & $C_5$ & $\{1\}$ & $0$ & $8$ & $4$\\
 &  &  & $\{1,9\}$ & $-2$ & $6$ & $1$\\
 &  &  & $(\Z/10\Z)^\times$ & $-3$ & $5$ & $-1/2$ \\
%
%
%
%
$G(1,12)$ & $12$ & $C_2$ & $\{1\}$ & $2$ & $10$ & $7$\\
 &  &  & $\{1,5\}$ & $3/2$ & $15/2$ & $21/4$\\
 &  &  & $\{1,7\}$ & $0$ & $8$ & $4$\\
 &  &  & $\{1,11\}$ & $3/2$ & $15/2$ & $21/4$\\
 &  &  & $(\Z/12\Z)^\times$ & $1/2$ & $13/2$ & $15/4$\\
& & $C_3$  & $\{1\}$ & $0$ & $10$ & $5$\\
 &  &  & $\{1,5\}$ & $-5/2$ & $15/2$ & $5/4$\\
 &  &  & $\{1,7\}$ & $0$ & $8$ & $4$\\
 &  &  & $\{1,11\}$ & $-3/2$ & $13/2$ & $7/4$ \\
 &  &  & $(\Z/12\Z)^\times$ & $-2$ & $6$ & $1$\\
%
%
%
%
$G(2,2)$ & $2$ & $C_2$ & $\{1\}$ & $-3/2$ & $5/2$ & $-1/4$\\
%
%
%
%
    $G(2,4)$ & $4$ &  $C_2$ & $\{1\}$ & $0$ & $4$ & $2$ \\
    &   & & $(\Z/4\Z)^\times$ & $-1$ & $3$ & $1/2$\\
%
%
%
%
    $G(2,6)$ & $6$ & $C_2$ &   $\{1\}$ & $-2$ & $6$ & $1$ \\
    &  &  & $(\Z/6\Z)^\times$ & $-5/2$ & $11/2$ & $1/4$\\
     &  & $C_3$ &  $\{1\}$ & $-1$ & $5$ & $3/2$ \\
    &  &  & $(\Z/6\Z)^\times$ & $-2$ & $4$ & $0$ \\ 
%
%
%
%
    $G(2,8)$ & $8$ & $C_2$ &   $\{1\}$ & $2$ & $10$ & $7$\\
    &  &  & $\{1, 3\}$ & $-1$ & $7$ & $5/2$\\
    &  &  & $\{1, 5\}$ & $0$ & $8$ & $4$\\
    &  &  & $\{1, 7\}$ & $-1$ & $7$ & $5/2$\\
    &  &  & $(\Z/8\Z)^\times$ & $-2$ & $6$ & $1$\\
%
%
%
%
    $G(3,3)$ & $3$ & $C_3$ & $\{1\}$ & $-2$ & $4$ & $0$\\
    &  &  & $(\Z/3\Z)^\times$ & $-3/2$ & $5/2$ & $-1/4$\\
%
%
%
%
    $G(3,6)$ & $6$ & $C_2$ & $\{1\}$ & $0$ & $8$ & $4$\\
    &  & & $(\Z/6\Z)^\times$ & $0$ & $6$ & $3$\\
    &  & $C_3$ & $\{1\}$ & $-4$ & $8$ & $0$\\
    &  & & $(\Z/6\Z)^\times$ & $-3$ & $5$ & $-1/2$\\
%
%
%
%
    $G(4,4)$ & $2$ & $C_2$ & $\{1\}$ & $-2$ & $4$ & $1$ \\
    & & & $(\Z/4\Z)^\times$ & $-2$ & $4$ & $0$\\
%
%
%
%
    $G(5,5)$ & $5$ & $C_5$ & $\{1\}$ & $-8$ & $12$ & $-2$ \\
    & &  & $\{1,4\}$ & $-5$ & $7$ & $-3/2$\\
     & &  & $(\Z/5\Z)^\times$ & $-7/2$ & $9/2$ & $-5/4$ \\
\bottomrule[1.5pt]
\end{longtblr}

\bibliographystyle{alpha}
\bibliography{bibfile}

@article {BKLPR15,
    AUTHOR = {Bhargava, Manjul and Kane, Daniel M. and Lenstra, Jr., Hendrik
              W. and Poonen, Bjorn and Rains, Eric},
     TITLE = {Modeling the distribution of ranks, {S}elmer groups, and
              {S}hafarevich-{T}ate groups of elliptic curves},
   JOURNAL = {Camb. J. Math.},
  FJOURNAL = {Cambridge Journal of Mathematics},
    VOLUME = {3},
      YEAR = {2015},
    NUMBER = {3},
     PAGES = {275--321},
      ISSN = {2168-0930,2168-0949},
   MRCLASS = {11G05 (14G25)},
  MRNUMBER = {3393023},
MRREVIEWER = {John\ L.\ Boxall},
       DOI = {10.4310/CJM.2015.v3.n3.a1},
       URL = {https://doi-org.dartmouth.idm.oclc.org/10.4310/CJM.2015.v3.n3.a1},
}

@misc{BS13a,
      title={The average number of elements in the 4-Selmer groups of elliptic curves is 7}, 
      author={Manjul Bhargava and Arul Shankar},
      year={2013},
      eprint={1312.7333},
      archivePrefix={arXiv},
      primaryClass={math.NT},
      note={Preprint, \texttt{arXiv:1312.7333}},
}

@misc{BS13b,
      title={The average size of the 5-selmer group of elliptic curves is 6, and the average rank is less than 1}, 
      author={Manjul Bhargava and Arul Shankar},
      year={2013},
      eprint={1312.7859},
      archivePrefix={arXiv},
      primaryClass={math.NT},
      note={Preprint, \texttt{arXiv:1312.7859}},
}

@article {BS15a,
    AUTHOR = {Bhargava, Manjul and Shankar, Arul},
     TITLE = {Binary quartic forms having bounded invariants, and the
              boundedness of the average rank of elliptic curves},
   JOURNAL = {Ann. of Math. (2)},
  FJOURNAL = {Annals of Mathematics. Second Series},
    VOLUME = {181},
      YEAR = {2015},
    NUMBER = {1},
     PAGES = {191--242},
      ISSN = {0003-486X},
   MRCLASS = {11E20 (14G25)},
  MRNUMBER = {3272925},
MRREVIEWER = {John M. Voight},
       DOI = {10.4007/annals.2015.181.1.3},
       URL = {https://doi.org/10.4007/annals.2015.181.1.3},
}

@article {BS15b,
    AUTHOR = {Bhargava, Manjul and Shankar, Arul},
     TITLE = {Ternary cubic forms having bounded invariants, and the
              existence of a positive proportion of elliptic curves having
              rank 0},
   JOURNAL = {Ann. of Math. (2)},
  FJOURNAL = {Annals of Mathematics. Second Series},
    VOLUME = {181},
      YEAR = {2015},
    NUMBER = {2},
     PAGES = {587--621},
      ISSN = {0003-486X},
   MRCLASS = {11G05 (11E20 14G05)},
  MRNUMBER = {3275847},
MRREVIEWER = {Joseph H. Silverman},
       DOI = {10.4007/annals.2015.181.2.4},
       URL = {https://doi.org/10.4007/annals.2015.181.2.4},
}

@article {Bil69,
    AUTHOR = {Billingsley, Patrick},
     TITLE = {On the central limit theorem for the prime divisor functions},
   JOURNAL = {Amer. Math. Monthly},
  FJOURNAL = {American Mathematical Monthly},
    VOLUME = {76},
      YEAR = {1969},
     PAGES = {132--139},
      ISSN = {0002-9890,1930-0972},
   MRCLASS = {60.30},
  MRNUMBER = {242222},
MRREVIEWER = {J.\ Sethuraman},
       DOI = {10.2307/2317259},
       URL = {https://doi.org/10.2307/2317259},
}

@article {BN22,
    AUTHOR = {Bruin, Peter and Najman, Filip},
     TITLE = {Counting elliptic curves with prescribed level structures over
              number fields},
   JOURNAL = {J. Lond. Math. Soc. (2)},
  FJOURNAL = {Journal of the London Mathematical Society. Second Series},
    VOLUME = {105},
      YEAR = {2022},
    NUMBER = {4},
     PAGES = {2415--2435},
      ISSN = {0024-6107,1469-7750},
   MRCLASS = {11G05 (11G18 11G50 11N45 14D23 14G40)},
  MRNUMBER = {4440538},
MRREVIEWER = {Robert\ Juricevic},
       DOI = {10.1112/jlms.12564},
       URL = {https://doi.org/10.1112/jlms.12564},
}

@inproceedings {CHL19,
    AUTHOR = {Chan, Stephanie and Hanselman, Jeroen and Li, Wanlin},
     TITLE = {Ranks, 2-{S}elmer groups, and {T}amagawa numbers of elliptic
              curves with {$\mathbb{Z}/2\mathbb{Z}\times\mathbb{Z}/8\mathbb{Z}$}-torsion},
 BOOKTITLE = {Proceedings of the {T}hirteenth {A}lgorithmic {N}umber
              {T}heory {S}ymposium},
    SERIES = {Open Book Ser.},
    VOLUME = {2},
     PAGES = {173--189},
 PUBLISHER = {Math. Sci. Publ., Berkeley, CA},
      YEAR = {2019},
      ISBN = {978-1-935107-03-3; 978-1-935107-02-6},
   MRCLASS = {11G05 (11Y40)},
  MRNUMBER = {3952011},
MRREVIEWER = {John\ L.\ Boxall},
}

@misc{CV25,
      title={Selmer groups of families of elliptic curves with an $\ell$-isogeny}, 
      author={Stephanie Chan and Matteo Verzobio},
      year={2025},
      eprint={2508.21406},
      archivePrefix={arXiv},
      primaryClass={math.NT},
      url={https://arxiv.org/abs/2508.21406}, 
}

@incollection {DDT94,
    AUTHOR = {Darmon, Henri and Diamond, Fred and Taylor, Richard},
     TITLE = {Fermat's last theorem},
 BOOKTITLE = {Current developments in mathematics, 1995 ({C}ambridge, {MA})},
     PAGES = {1--154},
 PUBLISHER = {Int. Press, Cambridge, MA},
      YEAR = {1994},
      ISBN = {1-57146-029-2},
   MRCLASS = {11G18 (11D41 11F80 11G05)},
  MRNUMBER = {1474977},
MRREVIEWER = {M.\ Ram\ Murty},
}

@article {DD15,
    AUTHOR = {Dokchitser, Tim and Dokchitser, Vladimir},
     TITLE = {Local invariants of isogenous elliptic curves},
   JOURNAL = {Trans. Amer. Math. Soc.},
  FJOURNAL = {Transactions of the American Mathematical Society},
    VOLUME = {367},
      YEAR = {2015},
    NUMBER = {6},
     PAGES = {4339--4358},
      ISSN = {0002-9947,1088-6850},
   MRCLASS = {11G07 (11G05 11G40)},
  MRNUMBER = {3324930},
MRREVIEWER = {Mahesh\ Kakde},
       DOI = {10.1090/S0002-9947-2014-06271-5},
       URL = {https://doi.org/10.1090/S0002-9947-2014-06271-5},
}

@article {EK40,
    AUTHOR = {Erd\"os, P. and Kac, M.},
     TITLE = {The {G}aussian law of errors in the theory of additive number
              theoretic functions},
   JOURNAL = {Amer. J. Math.},
  FJOURNAL = {American Journal of Mathematics},
    VOLUME = {62},
      YEAR = {1940},
     PAGES = {738--742},
      ISSN = {0002-9327,1080-6377},
   MRCLASS = {10.0X},
  MRNUMBER = {2374},
MRREVIEWER = {E.\ R.\ van Kampen},
       DOI = {10.2307/2371483},
       URL = {https://doi-org.dartmouth.idm.oclc.org/10.2307/2371483},
}

@article {Fis03,
    AUTHOR = {Fisher, Tom A.},
     TITLE = {The {C}assels-{T}ate pairing and the {P}latonic solids},
   JOURNAL = {J. Number Theory},
  FJOURNAL = {Journal of Number Theory},
    VOLUME = {98},
      YEAR = {2003},
    NUMBER = {1},
     PAGES = {105--155},
      ISSN = {0022-314X,1096-1658},
   MRCLASS = {11G05},
  MRNUMBER = {1950441},
MRREVIEWER = {Ivica\ Gusi\'c},
       DOI = {10.1016/S0022-314X(02)00038-0},
       URL = {https://doi.org/10.1016/S0022-314X(02)00038-0},
}

@misc{HL25,
      title={Genus formulas for families of modular curves}, 
      author={Asimina S. Hamakiotes and Jun Bo Lau},
      year={2025},
      eprint={2501.10883},
      archivePrefix={arXiv},
      primaryClass={math.NT},
      url={https://arxiv.org/abs/2501.10883}, 
}

@article {KLO14,
    AUTHOR = {Klagsbrun, Zev and Lemke Oliver, Robert J.},
     TITLE = {The distribution of the {T}amagawa ratio in the family of
              elliptic curves with a two-torsion point},
   JOURNAL = {Res. Math. Sci.},
  FJOURNAL = {Research in the Mathematical Sciences},
    VOLUME = {1},
      YEAR = {2014},
     PAGES = {Art. 15, 10},
      ISSN = {2522-0144,2197-9847},
   MRCLASS = {11G05 (14H25 14H52)},
  MRNUMBER = {3375649},
MRREVIEWER = {Andrew\ Bremner},
       DOI = {10.1186/s40687-014-0015-4},
       URL = {https://doi-org.dartmouth.idm.oclc.org/10.1186/s40687-014-0015-4},
}

@article {Kub76,
    AUTHOR = {Kubert, Daniel Sion},
     TITLE = {Universal bounds on the torsion of elliptic curves},
   JOURNAL = {Proc. London Math. Soc. (3)},
  FJOURNAL = {Proceedings of the London Mathematical Society. Third Series},
    VOLUME = {33},
      YEAR = {1976},
    NUMBER = {2},
     PAGES = {193--237},
      ISSN = {0024-6115,1460-244X},
   MRCLASS = {10B15 (14G25)},
  MRNUMBER = {434947},
MRREVIEWER = {J.\ W. S. Cassels},
       DOI = {10.1112/plms/s3-33.2.193},
       URL = {https://doi-org.dartmouth.idm.oclc.org/10.1112/plms/s3-33.2.193},
}

@misc{lmfdb,
  shorthand    = {LMFDB},
  author       = {The {LMFDB Collaboration}},
  title        = {The {L}-functions and modular forms database},
  howpublished = {\url{https://www.lmfdb.org}},
  year         = {2025},
  note         = {[Online; accessed 8 August 2025]},
}

@misc{LTG24+,
      title={The asymptotic estimation of prime ideals in imaginary quadratic fields and Chebyshev's bias}, 
      author={Chen Lin and Chenhao Tang and Xuejun Guo},
      year={2024},
      eprint={2412.14792},
      archivePrefix={arXiv},
      primaryClass={math.NT},
      url={https://arxiv.org/abs/2412.14792}, 
}

@book {Nar04,
    AUTHOR = {Narkiewicz, W{\l}adys{\l}aw},
     TITLE = {Elementary and analytic theory of algebraic numbers},
    SERIES = {Springer Monographs in Mathematics},
   EDITION = {Third},
 PUBLISHER = {Springer-Verlag, Berlin},
      YEAR = {2004},
     PAGES = {xii+708},
      ISBN = {3-540-21902-1},
   MRCLASS = {11Rxx (11-01 11-02)},
  MRNUMBER = {2078267},
       DOI = {10.1007/978-3-662-07001-7},
       URL = {https://doi.org/10.1007/978-3-662-07001-7},
}

@misc{Phi24+,
      title={Points of bounded height in images of morphisms of weighted projective stacks with applications to counting elliptic curves}, 
      author={Tristan Phillips},
      year={2024},
      eprint={2201.10624},
      archivePrefix={arXiv},
      primaryClass={math.NT},
      url={https://arxiv.org/abs/2201.10624}, 
}

@article {Phi25,
    AUTHOR = {Phillips, Tristan},
     TITLE = {Average {A}nalytic {R}anks of {E}lliptic {C}urves over
              {N}umber {F}ields},
   JOURNAL = {Forum Math. Sigma},
  FJOURNAL = {Forum of Mathematics. Sigma},
    VOLUME = {13},
      YEAR = {2025},
     PAGES = {Paper No. e40},
      ISSN = {2050-5094},
   MRCLASS = {11G05 (11D45 11G07 11G35 11G50 14G05)},
       DOI = {10.1017/fms.2024.127},
       URL = {https://doi.org/10.1017/fms.2024.127},
}

@article {PR12,
    AUTHOR = {Poonen, Bjorn and Rains, Eric},
     TITLE = {Random maximal isotropic subspaces and {S}elmer groups},
   JOURNAL = {J. Amer. Math. Soc.},
  FJOURNAL = {Journal of the American Mathematical Society},
    VOLUME = {25},
      YEAR = {2012},
    NUMBER = {1},
     PAGES = {245--269},
      ISSN = {0894-0347,1088-6834},
   MRCLASS = {11G10 (11G05 11G30 14G25 14K15)},
  MRNUMBER = {2833483},
MRREVIEWER = {Philipp\ Habegger},
       DOI = {10.1090/S0894-0347-2011-00710-8},
       URL = {https://doi-org.dartmouth.idm.oclc.org/10.1090/S0894-0347-2011-00710-8},
}

@article {Sch96,
    AUTHOR = {Schaefer, Edward F.},
     TITLE = {Class groups and {S}elmer groups},
   JOURNAL = {J. Number Theory},
  FJOURNAL = {Journal of Number Theory},
    VOLUME = {56},
      YEAR = {1996},
    NUMBER = {1},
     PAGES = {79--114},
      ISSN = {0022-314X,1096-1658},
   MRCLASS = {11G10 (11G30 14K02 14K15)},
  MRNUMBER = {1370197},
MRREVIEWER = {Philippe\ Satg\'e},
       DOI = {10.1006/jnth.1996.0006},
       URL = {https://doi.org/10.1006/jnth.1996.0006},
}

@book {Sil09,
    AUTHOR = {Silverman, Joseph H.},
     TITLE = {The arithmetic of elliptic curves},
    SERIES = {Graduate Texts in Mathematics},
    VOLUME = {106},
   EDITION = {Second},
 PUBLISHER = {Springer, Dordrecht},
      YEAR = {2009},
     PAGES = {xx+513},
      ISBN = {978-0-387-09493-9},
   MRCLASS = {11-02 (11G05 11G20 14H52 14K15)},
  MRNUMBER = {2514094},
MRREVIEWER = {Vasil\cprime\ \=I.\ Andr\=\i\u ichuk},
       DOI = {10.1007/978-0-387-09494-6},
       URL = {https://doi.org/10.1007/978-0-387-09494-6},
}

@inproceedings {Tat75,
    AUTHOR = {Tate, J.},
     TITLE = {Algorithm for determining the type of a singular fiber in an
              elliptic pencil},
 BOOKTITLE = {Modular functions of one variable, {IV} ({P}roc. {I}nternat.
              {S}ummer {S}chool, {U}niv. {A}ntwerp, {A}ntwerp, 1972)},
     PAGES = {33--52. Lecture Notes in Math., Vol. 476},
      YEAR = {1975},
   MRCLASS = {14G25 (14K15)},
  MRNUMBER = {0393039},
MRREVIEWER = {J. S. Milne},
}

@article {Vel71,
    AUTHOR = {Vélu, Jacques},
     TITLE = {Isog\'enies entre courbes elliptiques},
   JOURNAL = {C. R. Acad. Sci. Paris S\'er. A-B},
  FJOURNAL = {Comptes Rendus Hebdomadaires des S\'eances de l'Acad\'emie des
              Sciences. S\'eries A et B},
    VOLUME = {273},
      YEAR = {1971},
     PAGES = {A238--A241},
      ISSN = {0151-0509},
   MRCLASS = {14H99},
  MRNUMBER = {294345},
MRREVIEWER = {J.\ C.\ Fogarty},
}

\end{document}